\documentclass{article}
\usepackage{graphicx,pstricks}
\usepackage[all]{xy}
\usepackage{subfigure}
\usepackage{rotating}
\usepackage{amsfonts}
\usepackage{amssymb}
\usepackage{amsmath}
\usepackage{epsfig}
\usepackage{psfrag}
\usepackage{enumerate}
\usepackage{array}
\usepackage{theorem}

\newcommand{\C}{\mathbb{C}}
\newcommand{\R}{\mathbb{R}}
\newcommand{\N}{\mathbb{N}}
\newcommand{\Z}{\mathbb{Z}}

\newcommand{\Q}{\mathbb{Q}}

\newcommand{\Tr}{{\rm Tr}}
\newcommand{\hnc}{{\rm H}^n_\C}
\newcommand{\hncb}{\overline{{\rm H}}^n_\C}

\newtheorem{conj}{Conjecture}[section]
\newtheorem{thm}{Theorem}[section]

\newtheorem{lem}{Lemma}[section]
\newtheorem{prop}{Proposition}[section]
\newtheorem{dfn}{Definition}[section]

{\theorembodyfont{\rmfamily} \newtheorem{rmk}[thm]{Remark}}

\newcommand{\Pf}{{\em Proof}. }
\newcommand{\EPf}{\hfill$\Box$\vspace{.5cm}}
\numberwithin{equation}{section}

\date{December 13, 2010}

\title{Census of the complex hyperbolic sporadic triangle groups}
\author{Martin Deraux\\
Universit\'e de Grenoble I\\
Institut Fourier, B.P.74\\
38402 Saint-Martin-d'H\`eres Cedex\\
France \\
e-mail: {\tt deraux@ujf-grenoble.fr}
\and John R. Parker\\
Department of Mathematical Sciences\\
Durham University\\
South Road,\\
Durham DH1 3LE, England. \\
e-mail: {\tt j.r.parker@durham.ac.uk}
\and Julien Paupert\\
Department of Mathematics\\
University of Utah\\
155 South 1400 East\\
Salt Lake City, Utah 84112, USA.\\
e-mail: {\tt paupert@math.utah.edu}
}
        
\begin{document}
\maketitle

\begin{abstract}
  The goal of this paper is to give a conjectural census of complex
  hyperbolic sporadic triangle groups. We prove that only finitely many of
  these sporadic groups are lattices.

  We also give a conjectural list of all lattices among sporadic
  groups, and for each group in the list we give a conjectural group
  presentation, as well as a list of cusps and generators for their
  stabilisers. We describe strong evidence for these conjectural
  statements, showing that their validity depends on the solution of
  reasonably small systems of quadratic inequalities in four
  variables.
\end{abstract}

\section{Introduction}

The motivation for this paper is to construct discrete groups acting
on the complex hyperbolic plane ${\rm H}^2_\C$, more specifically
lattices (where one requires in addition that the quotient by the
action of the discrete group have finite volume). Complex hyperbolic
spaces ${\rm H}^n_{\mathbb{C}}$ are a natural generalisation to the
realm of K\"ahler geometry of the familiar non-Euclidean geometry of
${\rm H}^n_{\mathbb{R}}$. ${\rm H}^n_{\mathbb{C}}$ is simply the unit
ball in $\C^n$, endowed with the unique K\"ahler metric invariant
under all biholomorphisms of the ball; this metric is symmetric and
has non-constant negative real sectional curvature (holomorphic
sectional curvature is constant). The group of holomorphic isometries
of ${\rm H}^n_{\mathbb{C}}$ is the projectivised group ${\rm PU}(n,1)$
of a Hermitian form of Lorentzian signature $(n,1)$.

It is a well known fact due to Borel that lattices exist in the
isometry group of any symmetric space, but the general structure of
lattices and the detailed study of their representation theory brings
forth several open questions. The basic construction of lattices
relies on the fact that for any linear algebraic group $G$ defined
over $\mathbb{Q}$, the group of integral matrices $G(\mathbb{Z})$ is a
lattice in $G(\mathbb{R})$. $G(\mathbb{Z})$ is clearly discrete, and
the fact that it is a lattice follows from a theorem of Borel and
Harish-Chandra. More generally, to a group defined over a number field
(i.e. a finite extension of the rationals), one can associate a group
defined over $\mathbb{Q}$ by a process called restriction of
scalars. One is naturally led to the general notion of
\emph{arithmetic group}, keeping in mind that one would like to push
as far as possible the idea of taking integral matrices in a group
defined over $\mathbb{Q}$. For the general definition of
arithmeticity, we refer the reader to section~\ref{sec-arithm}. In the
context of the present paper, the arithmeticity criterion in that
section (Proposition~\ref{crit}) will be sufficient.

It is known since deep work of Margulis that lattices in the isometry
group of any symmetric space of higher rank (i.e. rank $\geqslant 2$) are
all arithmetic. There are four families of rank $1$ symmetric spaces
of non-compact type, namely
$$
{\rm H}^n_{\mathbb{R}}, {\rm H}^n_{\mathbb{C}}, {\rm H}^n_{\mathbb{H}}, {\rm H}^2_{\mathbb{O}}.
$$ 
Lattices in the isometry groups of the last two families (hyperbolic
spaces over the quaternions and the octonions) are all known to be
arithmetic, thanks to work of Corlette and Gromov-Schoen.

On the other hand, non-arithmetic lattices are known to exist in
${\rm PO}(n,1)$ (which is the isometry group of ${\rm H}^n_\mathbb{R}$) for
arbitrary $n\geqslant 2$. A handful of examples coming from Coxeter groups
were known in low dimensions before Gromov and Piatetski-Shapiro found
a general construction using so-called interbreeding of well-chosen
arithmetic real hyperbolic lattices (see~\cite{gps}).

The existence of non-arithmetic lattices in ${\rm PU}(n,1)$ (the group
of holomorphic isometries of ${\rm H}^n_\mathbb{C}$) for arbitrary $n$
is a longstanding open question. Examples are known only for
$n\leqslant 3$, and they are all commensurable to complex reflection
groups. More specifically, it turns out that all known non-arithmetic
lattices in ${\rm PU}(n,1)$ for $n=2$ or $3$ are commensurable to one
of the hypergeometric monodromy groups listed in~\cite{DM} and
\cite{M2} (the same list appears in~\cite{T}).

The goal of this paper is to announce (and give outstanding evidence
for) results that exhibit several new commensurability classes of
non-arithmetic lattices in ${\rm PU}(2,1)$. Our starting point was the
investigation by Parker and Paupert~\cite{vol2} of \emph{symmetric
  triangle groups}, i.e. groups generated by three complex reflections
of order $p\geqslant 3$ in a symmetric configuration (the case $p=2$
was studied by Parker in \cite{vol1}).

Writing $R_i$, $i=1,2,3$ for the generators, the symmetry condition means that there exists an isometry $J$ of
order $3$ such that $JR_iJ^{-1}=R_ {i+1}$ (indices mod $3$). It turns
out that conjugacy classes of symmetric triangle groups (with
generators of any fixed order $p\geqslant 2$) can then be parametrised by
$$
\tau={\rm Tr}(R_1J),
$$ 
provided we represent isometries by matrices for $R_1$ and $J$ in
${\rm SU}(2,1)$ (see section~\ref{sec-HnC} for basic geometric facts about
complex hyperbolic spaces).

Following~\cite{vol2}, we denote by $\Gamma(\frac{2\pi}{p},\tau)$ the
group generated by $R_1$ and $J$ as above. The main problem is to
determine the values $(p,\tau)$ of the parameters such that
$\Gamma(\frac{2\pi}{p},\tau)$ is a lattice in ${\rm PU}(2,1)$. It is a
difficult problem to do this in all generality (see the discussion
in~\cite{M1},~\cite{De1} for instance).

To simplify matters, we shall concentrate on a slightly smaller class
of groups. The results in~\cite{vol2} give the list of all values of
$p,\tau$ such that $R_1R_2$ and $R_1J$ are either parabolic, or
elliptic of finite order. When this condition holds, we refer to such
a triangle group as \emph{doubly elliptic} (see
section~\ref{sec-sporadic}).

It turns out that the double ellipticity condition is independent of
$p$, and the values of $\tau$ that yield doubly elliptic triangle
groups come into two continuous $1$-parameter families, together with
$18$ isolated values of the parameter $\tau$.

The continuous families yield groups that are subgroups of so-called
Mostow groups, i.e. ones where the generating reflections satisfy the
braid relation
$$
R_iR_{i+1}R_i=R_{i+1}R_iR_{i+1}.
$$ 
In that case, the problem of determining which parameters yield a
lattice is completely solved (see~\cite{M1},~\cite{M3} for the first
family and \cite{vol2} for the second).

The isolated values of $\tau$ corresponding to doubly elliptic
triangle groups are called \emph{sporadic values}, and the
corresponding triangle groups are called \emph{sporadic triangle
  groups} (the list of sporadic values is given in
Table~\ref{tab:sporadic}, page~\pageref{tab:sporadic}). It has been
suspected since~\cite{vol1} and~\cite{vol2} that sporadic
groups may yield interesting lattices. 

In fact, the work in~\cite{vol3} shows that only one sporadic triangle
group is an \emph{arithmetic} lattice; moreover, most sporadic
triangle groups are not commensurable to any of the previously known
non-arithmetic lattices (the Picard, Mostow and Deligne-Mostow
lattices). The precise statement of what ``most sporadic groups''
means is given in Theorem~\ref{comm}, see also~\cite{vol3}. The
question left open is of course to determine which sporadic groups are
indeed lattices.

To that end, it is quite natural to use the first author's computer
program (see~\cite{De1}), and to go through an experimental
investigation of the Dirichlet domains for sporadic groups. The goal
of the present paper is to report on the results of this search, which
turn out to be quite satisfactory.

We summarise the results of our computer experimentation in the
following (see section~\ref{sec-sporadic}, Table~\ref{tab:sporadic}
for the meaning of the parameters $\sigma_1,\dots,\sigma_9$):
\begin{conj}\label{conj:main}
  The following sporadic groups are non-arithmetic lattices in 
${\rm SU}(2,1)$:
\begin{itemize}
\item (cocompact): ${\bf \Gamma(\frac{2\pi}{5},\overline{\sigma}_4)}$,
  ${\bf \Gamma(\frac{2\pi}{8},\overline{\sigma}_4)}$, ${\bf
    \Gamma(\frac{2\pi}{12},\overline{\sigma}_4)}$.
 
\item (non cocompact): ${\bf \Gamma(\frac{2\pi}{3},\sigma_1)}$, ${\bf
    \Gamma(\frac{2\pi}{3},\sigma_5)}$,
  $\Gamma(\frac{2\pi}{4},\sigma_1)$, ${\bf
    \Gamma(\frac{2\pi}{4},\overline{\sigma}_4)}$, ${\bf
    \Gamma(\frac{2\pi}{4},\sigma_5)}$,
  $\Gamma(\frac{2\pi}{6},\sigma_1)$, ${\bf
    \Gamma(\frac{2\pi}{6},\overline{\sigma}_4)}$.
\end{itemize} 
\end{conj}

In fact we have obtained outstanding evidence that
Conjecture~\ref{conj:main} is correct, but this evidence was obtained
by doing numerical computations using floating point arithmetic, and
it is conceivable (though very unlikely) that the results are flawed
because of issues of precision, in a similar vein as the analysis
in~\cite{De1} of the results in~\cite{M1}. Instead of arguing that the
computer experimentation is not misleading, we will prove
Conjecture~\ref{conj:main} in \cite{vol5} by using more direct
geometric methods.

Note that the only part of Conjecture~\ref{conj:main} that is
conjectural is the fact that the groups in question are lattices. The
fact that these groups are not arithmetic follows from the results
in~\cite{vol2} and~\cite{vol3}. The groups indicated in bold are known
to not be commensurable to Deligne-Mostow-Picard lattices by
\cite{vol3} (in fact, for $\Gamma(\frac{2\pi}{4},\overline{\sigma}_4)$
and $\Gamma(\frac{2\pi}{6},\overline{\sigma}_4)$ this follows from
non-cocompactness by the arguments in \cite{vol3}).

Computer experiments also suggest that Conjecture~\ref{conj:main} is
essentially optimal. More specifically, sporadic groups that do not
appear in the list seem not to be lattices (most of them are not
discrete, a handful seem to have infinite covolume), apart from the
following:
\begin{equation}\label{eq:exceptions}
\Gamma\left(\frac{2\pi}{3},\overline{\sigma}_4\right),\ 
\Gamma\left(\frac{2\pi}{2},\sigma_5\right),\ 
\Gamma\left(\frac{2\pi}{2},\overline{\sigma}_5\right).\ 
\end{equation}
These exceptions are in fact completely understood, and they are all
arithmetic; the last two groups are both isomorphic to the lattice
studied in~\cite{De1} (see~\cite{vol1}). As for the first group,
partly thanks to work in~\cite{vol2}, we have:
\begin{thm} $\Gamma(\frac{2\pi}{3},\overline{\sigma}_4)$ is a
  cocompact arithmetic lattice in ${\rm SU}(2,1)$.
\end{thm}
The fact that this group is discrete was proved in~\cite{vol2}
(Proposition~6.4), the point being that all non-trivial Galois
conjugates of the relevant Hermitian form are definite. In fact it is
the only sporadic group that is contained in an arithmetic lattice,
by~\cite{vol3}. In order to check that it is cocompact, one uses the
same argument as in~\cite{De2}. More specifically, one needs to verify
that the Dirichlet domain is cocompact. This can be done without
knowing the precise combinatorics of that polyhedron (it is enough to
study a partial Dirichlet domain, and to verify that all the $2$-faces
of that polyhedron are compact, see \cite{De2}).

The non-discreteness results we prove in section~\ref{sec-nondiscrete}
of the paper are close to proving optimality of the statement of the
Conjecture, but the precise statement is somewhat lengthy (see
Theorem~\ref{thm:nondiscrete}). For now we simply state the following:
\begin{thm}
  Only finitely many sporadic triangle groups are discrete.
\end{thm}

\noindent
\textbf{Acknowledgements:} This project was funded in part by the NSF
grant DMS-0600816, through the funding of the first author's stay at
the University of Utah in September 2009. It is a pleasure to thank
Domingo Toledo for his interest and enthusiasm for this project.

\section{Arithmetic lattices arising from Hermitian forms over number
  fields}\label{sec-arithm}

For the sake of completeness, we recall in Def.~\ref{arithgen}
the general definition of arithmeticity (see also~\cite{zimmer},
chapter 6). For the purposes of the present paper the special case of
arithmetic groups arising from Hermitian forms over number fields will
be sufficient (see Proposition~\ref{crit} below).

Borel and Harish-Chandra proved that if $G$ is a linear algebraic
group defined over $\Q$ then $G(\Z)$ is a lattice in $G(\R)$. Recall
that a {\bf real linear algebraic group defined over $\Q$} is a
subgroup $G$ of ${\rm GL}(n,\R)$ for some $n$, such that the elements
of $G$ are precisely the solutions of a set of polynomial equations in
the entries of the matrices, with the coefficients of the polynomials
lying in $\Q$; one denotes $G(\R)=G$ and $G(\Z)=G \cap {\rm
  GL}(n,\Z)$.  From their result, one can deduce that any real
semisimple Lie group contains infinitely many (distinct
commensurability classes of) lattices, either cocompact or non
cocompact.  

One obtains the general definition by extending this notion to
all groups equivalent to groups of the form $G(\Z)$ in the following
sense:
\begin{dfn}\label{arithgen} Let $G$ be a semisimple Lie group, and
  $\Gamma$ a subgroup of $G$. Then $\Gamma$ is an \emph{arithmetic
    lattice} in $G$ if there exist an algebraic group $S$ defined over
  $\Q$ and a continuous homomorphism $\phi: S(\R)^0 \longrightarrow G$
  with compact kernel such that $\Gamma$ is commensurable to
  $\phi(S(\Z) \cap S(\R)^0)$.
\end{dfn}
The fact that $\Gamma$ as in the definition is indeed a
lattice follows from the Borel--Harish-Chandra theorem.

Here we focus on the case of \emph{integral groups arising from
  Hermitian forms over number fields}. This means that we consider
groups $\Gamma$ which are contained in ${\rm SU}(H,\mathcal{O}_K)$,
where $K$ is a number field, $\mathcal{O}_K$ denotes its ring of
algebraic integers, and $H$ is a Hermitian form of signature $(2,1)$
with coefficients in $K$. Note that $\mathcal{O}_K$ is usually not
discrete in $\mathbb{C}$, so $SU(H,\mathcal{O}_K)$ is usually not
discrete in $SU(H)$. Under an additional assumption on the Galois
conjugates $^\varphi H$ of the form (obtained by applying field
automorphisms $\varphi\in {\rm Gal}(K)$ to the entries of the
representative matrix of $H$), the group $SU(H,\mathcal{O}_K)$ is
indeed discrete (see part 1 of Prop.~\ref{crit}).

\begin{prop}\label{crit}
  Let $E$ be a purely imaginary quadratic extension of a totally real
  field $F$, and $H$ a Hermitian form of signature (2,1) defined over
  $E$. 
  \begin{enumerate}
  \item ${\rm SU}(H;\mathcal{O}_E)$ is a lattice in $SU(H)$ if and
    only if for all $\varphi \in {\rm Gal}(F)$ not inducing the
    identity on $F$, the form $^\varphi H$ is definite. Moreover, in
    that case, ${\rm SU}(H;\mathcal{O}_E)$ is an arithmetic lattice.
  \item Suppose $\Gamma\subset {\rm SU}(H;\mathcal{O}_E)$ is a
    lattice. Then $\Gamma$ is arithmetic if and only if $\varphi \in
    {\rm Gal}(F)$ not inducing the identity on $F$, the form $^\varphi
    H$ is definite.
  \end{enumerate}
\end{prop}  

Part 1 of the Proposition is quite natural (and motivates the
formulation of the general definition of arithmeticity). Indeed, it is
a general fact that one can embed $\mathcal{O}_K$ discretely into
$\C^r$ by 
$$
x\mapsto ({\varphi_1}(x),\dots,{\varphi_r}(x))
$$
where ${\varphi_1},\dots,{\varphi_r}$ denote the distinct embeddings of
$K$ into the complex numbers (up to complex conjugation).

The group $S=\prod_{j=1}^r SU(^{\varphi_j}H)$ can be checked to be
defined over $\Q$ (this is an instance of a general process called
\emph{restriction of scalars}). Its integer points correspond to
$\prod_{j=1}^r SU(^{\varphi_j}H,\mathcal{O}_K)$, which is a lattice in
$S(\R)$ by the theorem of Borel and Harish-Chandra.

Now the key point is that the assumption on the Galois conjugates
amounts to saying that the projection
$$
\prod_{j=1}^r SU(^{\varphi_j}H)\rightarrow SU(^{\varphi_1}H)
$$
onto the first factor has compact kernel, hence maps discrete sets to
discrete sets (compare with Definition~\ref{arithgen}). This implies
that $SU(H,\mathcal{O}_K)$ is a lattice in $SU(H)$.

The proof of part~2 of Proposition~\ref{crit} is a bit more
sophisticated (see lemma 4.1 of \cite{M1}, 12.2.6 of \cite{DM} or
Prop. 4.1 of \cite{vol3}). Note that when the group $\Gamma$ as
in the Proposition is non-arithmetic, it necessarily has infinite
index in $SU(H,\mathcal{O}_K)$ (which is non-discrete in $SU(H)$).

\section{Complex hyperbolic space and its isometries}\label{sec-HnC}

For the reader's convenience we include a brief summary of key
definitions and facts about complex hyperbolic geometry, see~\cite{G}
for more information.

Let $\langle \cdot,\cdot \rangle$ be a Hermitian form of signature
$(n,1)$ on $\C^{n+1}$, which we can describe in matrix form as
$$
\langle v,w\rangle=w^*Hv.
$$
The unitary group ${\rm U}(H)$ is the group of matrices that preserve this
inner product, i.e.
$$
{\rm U}(H)=\{M\in {\rm GL}(n+1,\C):M^*HM=H\}.
$$

The signature condition amounts to saying that after an appropriate
linear change of coordinates, the Hermitian inner product is the
standard Lorentzian Hermitian product
\begin{equation}\label{eq:stdform}
-v_0\overline{w}_0+v_1\overline{w}_1+\dots+v_n\overline{w}_n,
\end{equation}
whose unitary group is usually denoted by ${\rm U}(n,1)$.  For
computational purposes, it can be convenient to work with a
non-diagonal matrix $H$ (as we do throughout this paper), but of
course, under the $(n,1)$ signature assumption, ${\rm U}(H)$ is
isomorphic to ${\rm U}(n,1)$.

As a set, $\hnc$ is just the subset of projective space $P^n_\C$
corresponding to the set of negative lines in $\C^{n+1}$,
i.e. $\C$-lines spanned by a vector $v\in\C^{n+1}$ such that $\langle
v,v\rangle <0$. Working in coordinates where the form is diagonal, any
negative line is spanned by a unique vector of the form
$(1,v_1,\dots,v_n)$, and negativity translates into
$$
|v_1|^2+ \dots + |v_n|^2<1
$$
which shows how to describe complex hyperbolic space as the {\bf unit ball}
in $\C^{n}$.

It is often useful to consider the {\bf boundary} of complex
hyperbolic space, denoted by $\partial \hnc$. This corresponds to
the set of null lines, i.e. $\C$-lines spanned by nonzero vectors
$v\in\C^{n+1}$ with $\langle v,v\rangle =0$. In terms of the ball
model alluded to in the previous paragraph, the boundary is of course
simply the unit sphere in $\C^{n+1}$.

The group ${\rm PU}(H)$ clearly acts by biholomorphisms on $\hnc$ (the
action is effective and transitive), and it turns out that ${\rm
  PU}(H)$ is actually the group of all biholomorphisms of complex
hyperbolic space. There is a unique K\"ahler metric on $\hnc$
invariant under the action of ${\rm PU}(H)$ (it can be described as
the Bergman metric of the ball). We will not need any explicit formula
for the metric, all we need is the formula for the distance between
two points (this will be enough for the purposes of the present
paper). Writing $X,Y$ for negative vectors in $\C^{n+1}$ and $x,y$ for
the corresponding $\C$-lines in $\hnc$, we have
\begin{equation}\label{dist}
\cosh ^2\left(\frac{\rho\bigl(x,y\bigr)}{2}\right) 
= \frac{|\langle X, Y \rangle|^2}{\langle X, X \rangle  \langle Y, Y \rangle}
\end{equation}
The factor $1/2$ inside the hyperbolic cosine is included for purposes
of normalisation only (it ensures that the holomorphic sectional
curvature of ${\rm H}^n_\C$ is $-1$, rather than just any negative
constant).

It is not hard to see that
$${\rm  Isom}({\rm H}_\C^n)={\rm PU}(n,1) \rtimes \Z/2$$
where the $\Z/2$ factor corresponds to complex conjugation (any
involutive antiholomorphic isometry would do).

The usual {\bf classification of isometries} of negatively curved metric
spaces, in terms of the analysis of the fixed points in
$$
\hncb=\hnc\cup\partial \hnc,
$$
is used throughout in the paper. Any nontrivial $g \in {\rm PU}(n,1)$
is of precisely one of the following types:
\begin{itemize}
\item \emph{elliptic}: $g$ has a fixed point in ${\rm H}_\C^n$;
\item \emph{parabolic}: $g$ has exactly one fixed point in $\hncb$,
  which lies in $\partial{\rm H}_\C^n$;
\item \emph{loxodromic}: $g$ has exactly two fixed points in $\hncb$, which lie in $\partial{\rm H}_\C^n$.
\end{itemize}
In the special case $n=2$, there is a simple formula involving the
trace of a representative $G\in {\rm SU}(2,1)$ of $g\in {\rm PU}(2,1)$ to
determine the type of the isometry $g$ (see~\cite{G}, p.204).

We will sometimes use a slightly finer classification for elliptic
isometries, calling an element {\bf regular elliptic} if any of its
representative has pairwise distinct eigenvalues.  The eigenvalues of
a matrix $A \in {\rm U}(n,1)$ representing an elliptic isometry $g$ all have
modulus one. Exactly one of these eigenvalues has a eigenvector $v$
with $\langle v,v\rangle < 0$ (the span of $v$ gives a fixed point of $g$
in ${\rm H}_\C^n$), and such an eigenvalue will be called \emph{of
  negative type}. Regular elliptic isometries have an isolated fixed
point in ${\rm H}_\C^n$.

Among non-regular elliptic elements, one finds {\bf complex
  reflections}, whose fixed point sets are totally geodesic copies of
${\rm H}^{n-1}_\C$ embedded in $\hnc$. More specifically, such
``complex hyperplanes'' can be described by a positive line in
$\C^{n+1}$, i.e. a $\C$-line spanned by a vector $v$ with $\langle
v,v\rangle >0$. Given such a vector, the set of $\C$-lines contained
in
$$
v^{\perp}=\left\{ w\in\C^{n+1}:\langle v,w\rangle =0 \right\}
$$
intersects $\hnc$ in a copy of ${\rm H}^{n-1}_\C$. The point in
projective space corresponding to $v$ is called \emph{polar} to the
hyperplane determined by $v^\perp$. In terms of the ball model, these
copies of ${\rm H}^{n-1}_\C$ simply correspond to the intersection
with the unit ball of affine hyperplanes in $\C^n$. If $v$ is a
positive vector, any isometry of $\hnc$ fixing the lines in $v^\perp$,
can be described as
$$
x\mapsto x+(\zeta-1)\frac{\langle x,v\rangle}{\langle v,v\rangle}v
$$
for some $\zeta\in\C$ with $|\zeta|=1$. The corresponding isometry is
called a complex reflection, $\zeta$ is called its multiplier,
and the argument of $\zeta$ is referred to as the rotation angle of
the complex reflection.

Note that the respective positions of two complex hyperplanes are
easily read off in terms of their polar vectors. Indeed, we have the
following (see~\cite{G}, p.100).
\begin{lem}\label{lem:twolines}
  Let $v_1,v_2$ be positive vectors in $\C^{n+1}$, and let $L_1,L_2$
  denote the corresponding complex hyperplanes in $\hnc$. Let 
  $$
  C = \frac{|\langle v_1,v_2\rangle|^2}{\langle v_1,v_1\rangle \langle v_2,v_2\rangle}.
  $$
\begin{itemize}
\item[(1)] $L_1$ and $L_2$ intersect in ${\rm H}_\C^n$ $\iff$
  $C<1$. In that case the angle $\theta$ between $L_1$ and $L_2$
  satisfies $\cos \theta = C$.
\item[(2)] $L_1$ and $L_2$ intersect in $\partial {\rm H}_\C^n$ $\iff$
  $C=1$.
\item[(3)] $L_1$ and $L_2$ are ultraparallel $\iff$ $C>1$. In that
  case the distance $\rho$ between $L_1$ and $L_2$ satisfies $\cosh
  \frac{\rho}{2} = C$.
\end{itemize}
\end{lem}
Lemma~\ref{lem:twolines} will be used to get the discreteness test in
section~\ref{sec-nondiscrete} (the complex hyperbolic J\o rgensen's
inequality established in \cite{JKP}).

Parabolic isometries are either {\bf unipotent} or {\bf screw
  parabolic}; in the former case they are also called {\bf Heisenberg
  translations} (because the group of unipotent isometries fixing a
given point on $\partial \hnc$ is isomorphic to the Heisenberg
group $\mathcal{H}^{2n-1}$). There are two conjugacy classes of Heisenberg
translations, the {\bf vertical translations} (corresponding to the
centre, or commutator subgroup, of the Heisenberg group) and the
{\bf non-vertical translations} (see~\cite{G} for more details on
this discussion).

\section{Sporadic groups} \label{sec-sporadic} 

In this section we setup some notation and recall the main results
from \cite{vol2} and \cite{vol3}.
\begin{dfn}
  A \emph{symmetric triangle group} is a group generated by two
  elements $R_1$, $J\in SU(2,1)$ where $R_1$ is a complex reflection
  of order $p$ and $J$ is a regular elliptic isometry $J$ of order
  $3$.
\end{dfn}
The reason we call this a triangle group is that it is a subgroup of
index at most three in the group generated by three complex
reflections $R_1$, $R_2$ and $R_3$, defined by
\begin{equation}\label{tgroup}
R_2=JR_1J^{-1},\quad R_3=JR_2J^{-1},
\end{equation}
and we think of their three mirrors as describing a ``triangle'' of
complex lines (however the mirrors of the various $R_j$ need not
intersect in general).

The basic observation is that symmetric triangle groups can be
parameterised up to conjugacy by the order of $p$ of $R_1$ and
$$
\tau=\Tr(R_1J).
$$
We denote by $\psi=2\pi/p$ the rotation angle of $R_1$, and by
$$
\Gamma(\psi,\tau)
$$ 
the group generated by a complex reflection $R_1,J$ as above.

The generators for this group can be described explicitly by matrices of
the form
\begin{eqnarray}
J & = & \left[\begin{matrix} 0 & 0 & 1 \\ 1 & 0 & 0 \\ 0 & 1 & 0 
\end{matrix}\right] \label{eq-J}\\
R_1 & = & 
\left[\begin{matrix} e^{2i\psi/3} & \tau & 
-e^{i\psi/3}\,\overline{\tau} \\
0 & e^{-i\psi/3} & 0 \\ 0 & 0 & e^{-i\psi/3} \end{matrix}\right]
\label{eq-R1} 
\end{eqnarray}
These preserve the Hermitian form
$\langle {\bf z},{\bf w}\rangle={\bf w}^*H_\tau{\bf z}$
where
\begin{equation}\label{eq-H-tau}
H_\tau=\left[\begin{matrix} 
2\sin(\psi/2) & -ie^{-i\psi/6}\tau & ie^{i\psi/6}\overline{\tau} \\
ie^{i\psi/6}\overline{\tau} & 2\sin(\psi/2) & -ie^{-i\psi/6}\tau \\
-ie^{-i\psi/6}\tau & ie^{i\psi/6}\overline{\tau} & 2\sin(\psi/2) 
\end{matrix}\right].
\end{equation}  
The above matrices always generate a subgroup $\Gamma$ of ${\rm
  GL}(3,\C)$, but the signature of $H_\tau$ depends on the values of
$\psi$ and $\tau$. For any fixed value of $\psi$, the parameter space
for $\tau$ is described in Sections~2.4 and~2.6 of~\cite{vol2}.
\begin{dfn}
  The symmetric triangle group generated by $R_1$ and $J$ as
  in~\eqref{eq-J} and~\eqref{eq-R1} is called \emph{hyperbolic} is
  $H_\tau$ has signature $(2,1)$.
\end{dfn}

In order to get a tractable class of groups, we shall assume that
$R_1J$ is elliptic, and that $R_1R_2=R_1JR_1J^{-1}$ is either elliptic
or parabolic. The motivation for this condition is explained
in~\cite{vol1},~\cite{vol2} (it is quite natural in the context of the
search for \emph{lattices}, rather than discrete groups of possibly
infinite covolume).

A basic necessary condition for a subgroup of ${\rm PU}(2,1)$ to be
discrete is that all its elliptic elements must have finite
order, hence we make the following definition.
\begin{dfn}
  A symmetric triangle group is called \emph{doubly elliptic} if
  $R_1J$ is elliptic of finite order and $R_1R_2=R_1JR_1J^{-1}$
  is either elliptic of finite order or parabolic.
\end{dfn}

The list of parameters that yield double elliptic triangle groups was
obtained in~\cite{vol1} (see also~\cite{vol2}), by using a result of
Conway and Jones on sums of roots of unity. We recall the result in
the following.
\begin{thm}\label{discnec}
  Let $\Gamma$ be a symmetric triangle group such that $R_1J$ is
  elliptic and $R_1R_2$ is either elliptic or parabolic. If $R_1J$ and
  $R_1R_2$ have finite order (or are parabolic), then one of the
  following holds:
\begin{itemize}
\item $\Gamma$ is one of Mostow's lattices ($\tau=e^{i\phi}$ for some
  $\phi$).
\item $\Gamma$ is a subgroup of one of Mostow's lattices
  ($\tau=e^{2i\phi}+e^{-i\phi}$ for some $\phi$).
\item $\Gamma$ is one of the sporadic triangle groups, i.e $\tau \in
  \{\sigma_1,\overline{\sigma}_1,...,\sigma_9,\overline{\sigma}_9 \}$
  where the $\sigma_j$ are given in Table~\ref{tab:sporadic}.
\end{itemize}
\end{thm} 

\begin{table}[htbp]
$\begin{array}{l|l|l}
\sigma_1=e^{i\pi/3}+e^{-i\pi/6}\,2\cos(\pi/4) 
&
\sigma_2=e^{i\pi/3}+e^{-i\pi/6}\,2\cos(\pi/5)
&
\sigma_3=e^{i\pi/3}+e^{-i\pi/6}\,2\cos(2\pi/5)\\

\sigma_4=e^{2\pi i/7}+e^{4\pi i/7}+e^{8\pi i/7}
&
\sigma_5=e^{2\pi i/9}+e^{-i\pi/9}\,2\cos(2\pi/5)
&
\sigma_6=e^{2\pi i/9}+e^{-\pi i/9}\,2\cos(4\pi/5)\\

\sigma_7=e^{2\pi i/9}+e^{-i\pi/9}\,2\cos(2\pi/7)
&
\sigma_8=e^{2\pi i/9}+e^{-i\pi/9}\,2\cos(4\pi/7)
&
\sigma_9=e^{2\pi i/9}+e^{-i\pi/9}\,2\cos(6\pi/7)
\end{array}$
\caption{The 18 sporadic values are given by $\sigma_j$ or
  $\overline{\sigma}_j$, $j=1,\dots,9$. They correspond to isolated
  values of the parameter $\tau$ for which any
  $\Gamma(\frac{2\pi}{p},\tau)$ is doubly elliptic,  i.e. $R_1R_2$ and
  $R_1J$ are either parabolic or elliptic of finite
  order.}\label{tab:sporadic}
\end{table}

Therefore, for each value of $p\geqslant 3$, we have a finite number
of new groups to study, the $\Gamma(2\pi/p,\sigma_i)$ and
$\Gamma(2\pi/p,\overline{\sigma_i})$ which are hyperbolic. 
The list of sporadic groups that are hyperbolic is
given in the table of Section~3.3 of~\cite{vol2} 
(and we give them below in Table \ref{tab:nondisc}); for the sake of
brevity we only recall the following:

\begin{prop} For $p\geqslant 4$ and
  $\tau=\sigma_1,\sigma_2,\sigma_3,\overline{\sigma}_4,\sigma_5,\sigma_6,
\sigma_7,\overline{\sigma}_8$  or $\sigma_9$, $\Gamma(2\pi/p,\tau)$ 
is hyperbolic.
\end{prop}  

It was shown in \cite{vol2} that some of the hyperbolic sporadic
groups are non-discrete (see Corollary 4.2, Proposition 4.5 and
Corollary 6.4 of \cite{vol2}), essentially by using the lists of
discrete triangle groups on the sphere, the Euclidean plane and the
hyperbolic plane (this list is due to Schwarz in the spherical case,
and to Knapp in the hyperbolic case). For the convenience of the
reader, we recall the main non-discreteness results from~\cite{vol2} in the following:
\begin{prop}\label{nondiscrete} For $p \geqslant 3$ and ($\tau$ or $
  \overline{\tau}=\sigma_3,\sigma_8$ or $\sigma_9$),
  $\Gamma(2\pi/p,\tau)$ is not discrete. Also, for $p \geqslant 3,\
  p\neq 5$ and ($\tau$ or $\overline{\tau}=\sigma_6$),
  $\Gamma(2\pi/p,\tau)$ is not discrete.
\end{prop}

The new non-discreteness results contained in
section~\ref{sec-nondiscrete} push the same idea much further, by a
series of technical algebraic manipulations (in some places we use J\o
rgensen's inequality and a complex hyperbolic version of Shimizu's
lemma due to the second author, see Theorem~\ref{thm:shim});

The main results of \cite{vol3} are the following two statements. The
first result was obtained by applying the arithmeticity criterion from
Proposition~\ref{crit}. The second result was obtained by finding a
commensurability invariant which distinguishes the various groups
$\Gamma$, namely the field $\Q [{\rm Tr Ad} \Gamma]$ (the trace field
of the adjoint representation of $\Gamma$).
\begin{thm}\label{nonar}  Let $p\geqslant 3$ and $\tau\in
  \{\sigma_1,\overline{\sigma}_1,...,\sigma_9,\overline{\sigma}_9 \}$,
  and suppose that the triangle group $\Gamma(2\pi/p, \tau)$ is
  hyperbolic, and that it is a lattice in $SU(H_\tau)$.  Then
  $\Gamma(2\pi/p, \tau)$ is arithmetic if and only if $p=3$ and
  $\tau=\overline{\sigma}_4$.
\end{thm}

\begin{thm}\label{comm} The sporadic groups $\Gamma(2\pi/p,\tau)$
  ($p\geqslant3$ and $\tau\in
  \{\sigma_1,\overline{\sigma}_1,...,\sigma_9,\overline{\sigma}_9\}$)
  fall into infinitely many distinct commensurability
  classes. Moreover, they are not commensurable to any Picard or
  Mostow lattice, except possibly when:
$$
\begin{array}{lll}
 \bullet \ p=4 \ {\rm or} \ 6 &
 \bullet \ p=3 \ {\rm and} \ \tau=\sigma_7 &
 \bullet \ p=5 \ {\rm and} \ \tau \ {\rm or} \ \overline{\tau}=\sigma_1,\sigma_2 \\
 \bullet \ p=7 \ {\rm and} \ \tau=\overline{\sigma}_4 &
 \bullet \ p=8 \ {\rm and} \ \tau=\sigma_1 &
 \bullet \ p=10 \ {\rm and} \ \tau=\sigma_1,\sigma_2,\overline{\sigma}_2 \\
 \bullet \ p=12 \ {\rm and} \ \tau=\sigma_1,\sigma_7 &
 \bullet \ p=20 \ {\rm and} \ \tau=\sigma_1,\sigma_2 &
 \bullet \ p=24 \ {\rm and} \ \tau=\sigma_1
\end{array}
$$
\end{thm}

\section{Dirichlet domains}

Given a subgroup $\Gamma$ of ${\rm PU}(2,1)$, the Dirichlet domain for
$\Gamma$ centred at $p_0$ is the set:
$$ 
F_\Gamma = \ \left\{ \, x\in {\rm H}_\C^2 : d(x,p_0) \leqslant
  d(x,\gamma p_0), \forall \gamma \in \Gamma \right\}.
$$
A basic fact is that $\Gamma$ is discrete if and only if $F_\Gamma $
has nonempty interior, and in that case $F_\Gamma$ is a fundamental
domain for $\Gamma$ modulo the action of the (finite) stabiliser of
$p_0$ in $\Gamma$.

The simplicity of this general notion, and its somewhat canonical
nature (it only depends on the choice of the centre $p_0$), make
Dirichlet domains convenient to use in computer investigation as in
\cite{M1}, \cite{Ri}, \cite{De1} and \cite{De2}. Note however that
there is no algorithm to decide whether the set $F_\Gamma$ has
non-empty interior, and the procedure we describe below may never end
(this is already the case in the constant curvature setting, i.e. in
real hyperbolic space of dimension at least 3, see for
instance~\cite{EP}).

Our computer search is quite a bit more delicate than the search for
fundamental domains in the setting of arithmetic groups. The recent
announcement that Cartwright and Steger have been able to find
presentations for the fundamental groups of all so-called \emph{fake
  projective planes} mentions the use of massive computer calculations
in the same vein as our work (see~\cite{CS}), but there are major
differences however.

They use Dirichlet domains, but their task is facilitated by the fact
that the fundamental groups of fake projective planes are known to be
\emph{arithmetic} subgroups of ${\rm PU}(2,1)$ (see~\cite{K}
and~\cite{Y}). In particular, all the groups they consider are known
to be discrete \emph{a priori} (which is certainly not the
case for most complex hyperbolic sporadic groups). Cartwright and
Steger also use the knowledge of the volumes of the corresponding
fundamental domains (the list of arithmetic lattices that could
possibly contain the fundamental group of a fake projective plane is
brought down to a finite list by using Prasad's volume formula
\cite{Pra}). This allows one to check whether a partial Dirichlet
domain
$$ 
F_W = \ \left\{ \, x\in {\rm H}_\C^2 : d(x,p_0) \leqslant
  d(x,\gamma p_0), \forall \gamma \in W \right\}
$$
determined by a given finite set $W\subset\Gamma$ is actually equal to
$F_\Gamma$.

For an arbitrary discrete subgroup $\Gamma\subset {\rm PU}(2,1)$ and an
arbitrary choice of the centre $p_0$, the set $F_\Gamma$ is a
polyhedron bounded by \emph{bisectors} (see~\cite{M1} and~\cite{G}),
but it may have infinitely many faces, even if $\Gamma$ is
geometrically finite (see~\cite{B}).

Moreover, the combinatorics of Dirichlet domains tend to be
unnecessarily complicated, and one usually expects that simpler
fundamental domains can be obtained by suitable clever geometric
constructions. This general idea is illustrated by Dirichlet domains
for lattices in $\R^2$: when the group is not a rectangular lattice,
i.e. not generated by two translations along orthogonal axes, the
Dirichlet domain centred at any point is a hexagon (rather than a
parallelogram).

In ${\rm H}_\C^2$, Dirichlet domains typically contain digons (pairs
of vertices connected by distinct edges), see Figure~\ref{fig:pics3s4c}. In
particular the 1-skeleton is not piecewise totally geodesic. One can
also check that the 2-faces of a Dirichlet domain can never be
contained in a totally real totally geodesic copy of ${\rm H}^2_\R$, which
makes this notion a little bit unnatural (this was part of the
motivation behind the constructions of~\cite{DFP}, where fundamental
domains with simpler combinatorics that those in~\cite{M1} were
obtained).

\section{Experimental results}

\subsection{The G-procedure} \label{sec-gprocedure}

In order to sift through the complex hyperbolic sporadic groups, we
have run the procedures explained in~\cite{De1} and~\cite{De2} in
order to explore the Dirichlet domains centred at the centre of mass
of the mirrors of the three generating reflexions.

In terms of the notation in Section~\ref{sec-sporadic}, we take $p_0$
to be the unique fixed point in ${\rm H}^2_\C$ of the regular elliptic
element $J$ (this point is given either by $(1,1,1)$,
$(1,\omega,\overline{\omega})$ or $(1,\overline{\omega},\omega)$ for
$\omega=(-1+i\sqrt{3})/2$, depending on the parameters $p$ and $\tau$).

We start with the generating set $W_0=\{R_1^{\pm 1},R_2^{\pm
  1},R_3^{\pm 1}\}$ for $\Gamma$, and construct an increasing sequence
of sets $W_0\subset W_1\subset W_2\subset\dots$ by the G-procedure
(named after G. Giraud, see~\cite{De1} for the explanation of this
terminology).

First define a \emph{G-step} of the procedure by:
$$
G(W) = W \cup \{\alpha^{-1}\beta : \alpha,\beta\in W \textrm{
  yield a non-empty generic $2$-face of $F_{W}$} \}
$$

Here ``yielding a non-empty $2$-face of $F_{W}$'' means that the set
of points of $F_{W}$ that are equidistant from $p_0$, $\alpha p_0$
and $\beta p_0$ has dimension two (i.e. it has non-empty interior in
the corresponding intersection of two bisectors). ``Generic'' means
that this $2$-face is not contained in a complex geodesic
(see~\cite{De1}).

\begin{dfn}
  The set $W$ is called G-closed if $G(W)=W$.
\end{dfn}

The sequence $W_k$ is defined inductively by 
$$W_{k+1}=G(W_k).$$
The hope is that this sequence stabilises to a G-closed set $W=W_N$
after a finite number of steps. In particular, this procedure is
probably suitable only for the search for lattices (not for discrete
groups with infinite covolume).

\subsection{Issues of precision} \label{sec-precision}

The determination of the sequence of sets $W_k$ described
in Section~\ref{sec-gprocedure} depends on being able to determine the precise
list of all nonempty $2$-faces of the polyhedron $W$, for a given
finite set $W\subset\Gamma$. The difficult part is to prove that two
bisectors really yield a subset of $F_W$ of dimension smaller than
$2$, when they appear to do so numerically. 

Recall that the polyhedron $F_W$ is described by a (possibly large)
set of quadratic inequalities in $4$ variables (the real and imaginary
parts of the ball coordinates, for instance), where the coefficients
of the quadratic polynomials are obtained from matrices which are
possibly very long words in the generators $R_1,R_2,R_3$.

The computation of these matrices can be done without loss of
precision, since it can be reduced to arithmetic in the relevant
number field (see Section~2.5 of~\cite{vol2}). 

It is not clear how to solve the corresponding system of quadratic
inequalities. In order to save computational time, and for the lack of
having better methods, we have chosen to do all the computations
numerically, with a fixed (somewhat rough) precision, essentially the
same way as described in~\cite{De2}. We now briefly summarise what our
computer program does.

For a given (coequidistant) bisector intersection $B$, we need a
method to test whether $B\cap F_W$ has dimension two. In order to do
this, we work in spinal coordinates (see~\cite{De1}), and fit the disk
$B$ into a rectangular $N\times N$ grid. The $2$-face is declared
non-empty whenever we find more than one point in a given horizontal
and in a given vertical line in the grid. For the default version of
the program, we take $N=1000$.

In particular, the above description suggests that whenever the
polyhedron $F_W$ becomes small enough, our program will not find any
$2$-face whatsoever. If this happens at some stage $k$, the program
will consider $W_k$ as being G-closed and stop.

When fed a group that has infinite covolume, one expects that the
program would often run forever, since in that case Dirichlet domains
have tend to infinitely many faces. In practice, after a certain
number of steps, the sets $W_k$ are too large for the computer's
capacity, and the program will crash.

For the groups we have tested (namely all sporadic groups with $p\leqslant
24$), we have found these three behaviours:
\begin{enumerate}[A:]
\item The program finds a G-closed set $W_N=G(W_N)$, and the set
of numerically non-empty two-faces is non-empty.
\item The program finds a set $W_N$ for which it does not find any
  nonempty $2$-face whatsoever (in particular $W_N$ is Giraud closed,
  so the program stops).
\item The program exceeds its capacity in memory and crashes.
\end{enumerate}

As a working hypothesis, we shall interpret Behaviour~B as meaning that
the group is not discrete, and Behaviour~C meaning that the group has
infinite covolume (the latter behaviour is of course also conceivable
when the group is actually not discrete, or when we make a bad choice
of the centre of the Dirichlet domain).

\subsection{Census of sporadic groups generated by reflections of small order}

The computer program available on the first author's webpage at
\begin{center}
  {\tt http://www-fourier.ujf-grenoble.fr/$\sim$deraux/java}
\end{center}
was run for all sporadic groups (see Section~\ref{sec-sporadic}) with
$2\leqslant p\leqslant 24$.

The groups with $p=2$ were analysed by Parker in~\cite{vol1}, and our
program confirms his results; in that case $\tau$ and
$\overline{\tau}$ give the same groups, and only $\tau=\sigma_5$ or
$\sigma_7$ appear to be discrete. Both exhibit Behaviour~A, but the
first one gives a compact polyhedron; as mentioned in the
introduction, this lattice is actually the same as the
$(4,4,4;5)$-triangle group, i.e. the group that is studied
in~\cite{De2}, see~\cite{vol1} and \cite{Sc}. The Giraud-closed
polyhedron obtained for $\sigma_7$ has infinite volume.

For $3\leqslant p\leqslant 24$, there are few groups that exhibit
Behaviour~A (as defined in Section~\ref{sec-precision}), namely: all
groups with $\tau=\overline{\sigma}_4$, those with $\tau=\sigma_1$,
$p=3,4,5,6$, and finally those with $\tau=\sigma_5$, $p=3,4$ or $5$.

\begin{figure}[htbp]
  \centering
  \hfill\subfigure[$1$]{
    \epsfig{figure=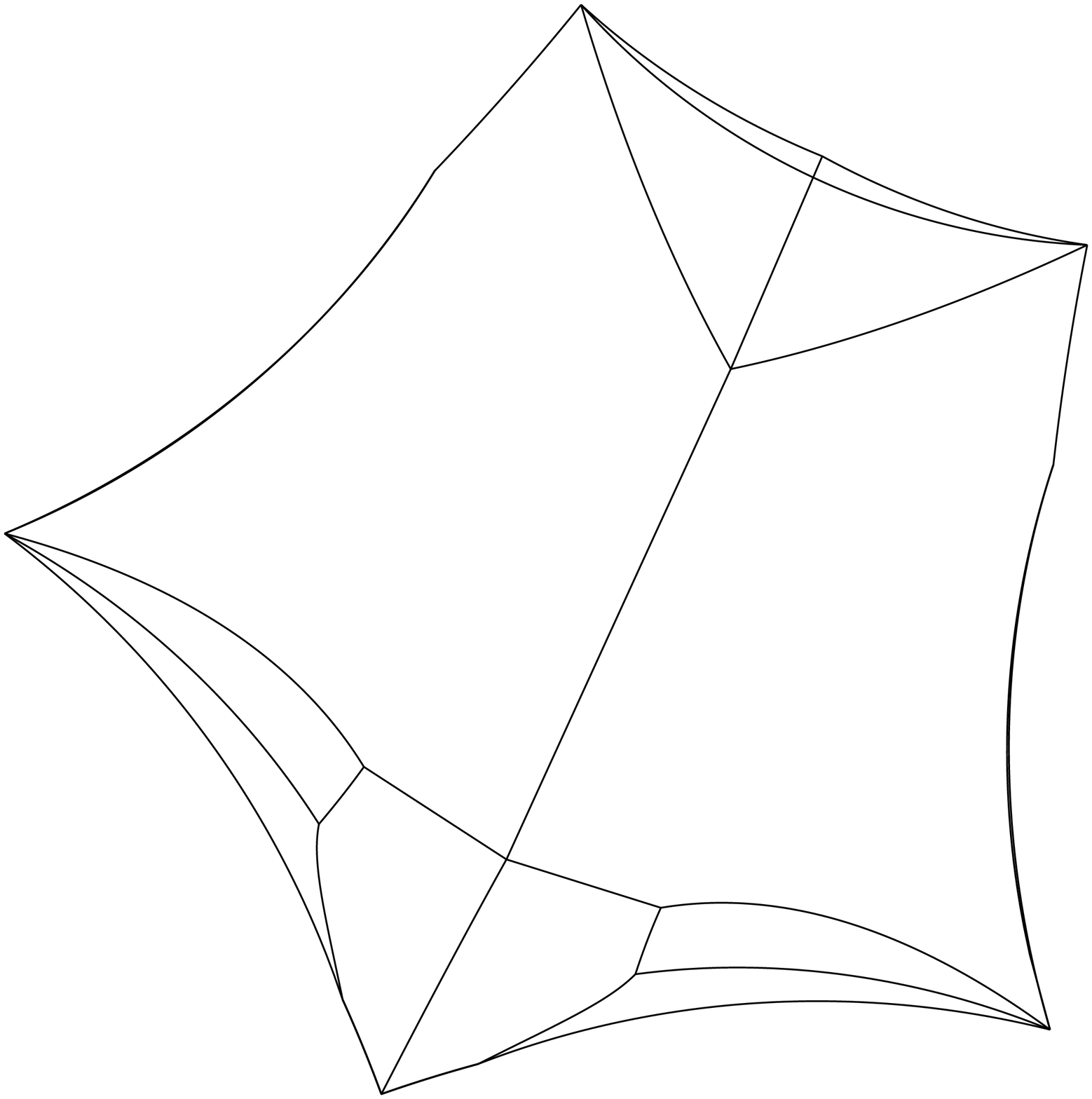, width=0.2\textwidth} }\hfill
  \subfigure[$12$]{
    \epsfig{figure=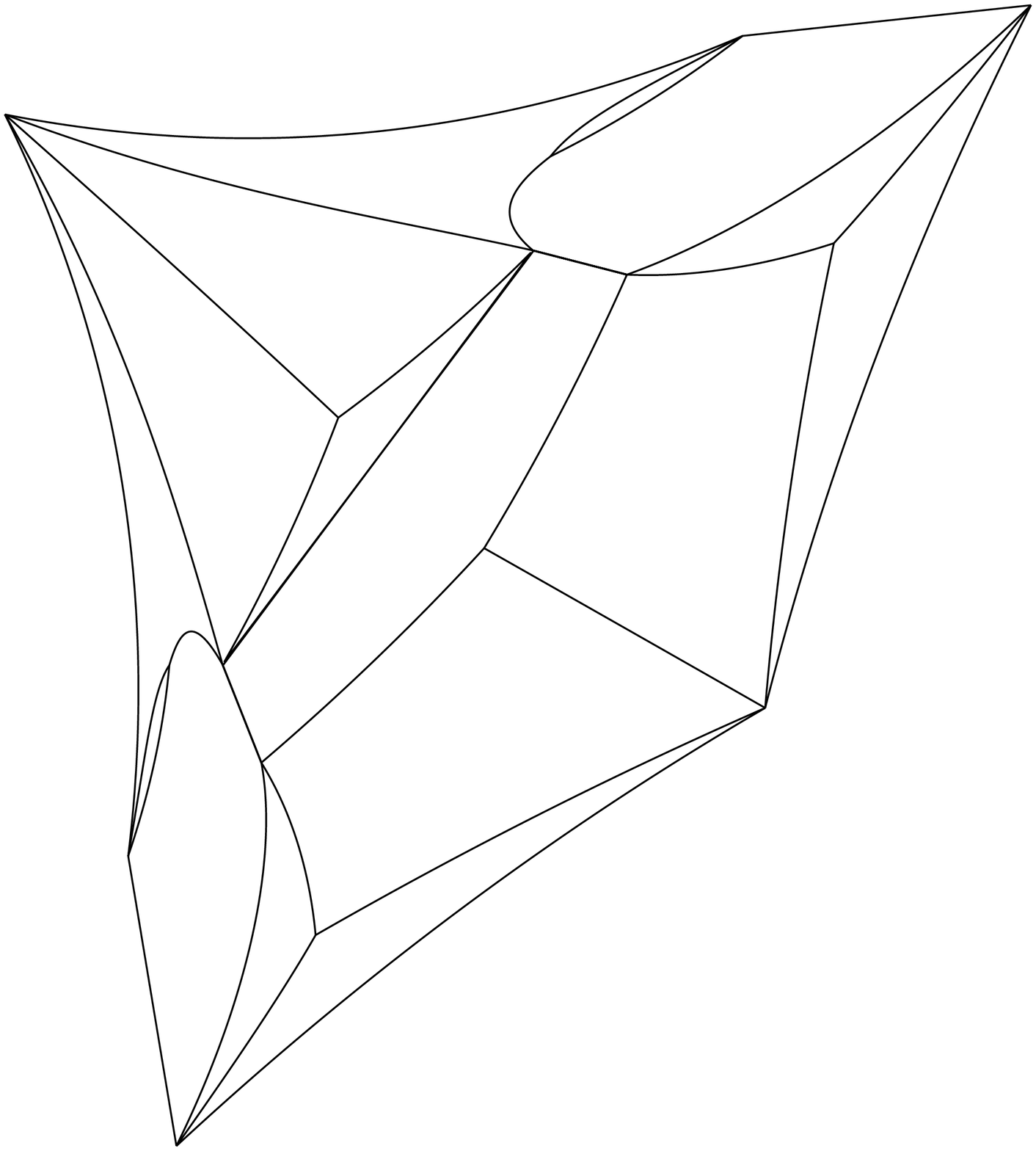, width=0.2\textwidth} }\hfill\mbox{}

  \hfill
  \subfigure[$121$]{
    \epsfig{figure=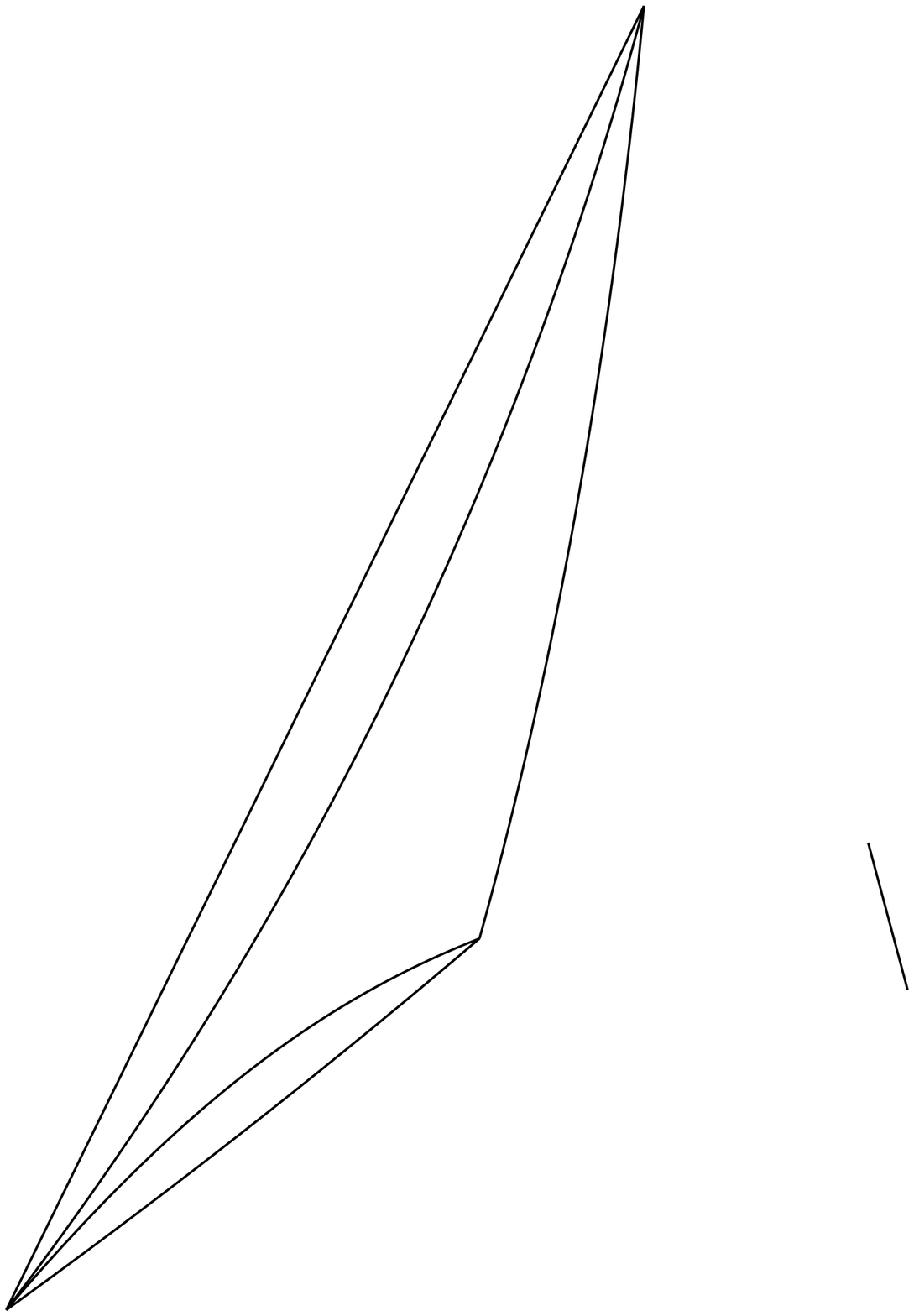, width=0.1\textwidth} }\hfill
  \subfigure[$12\overline{1}$]{
    \epsfig{figure=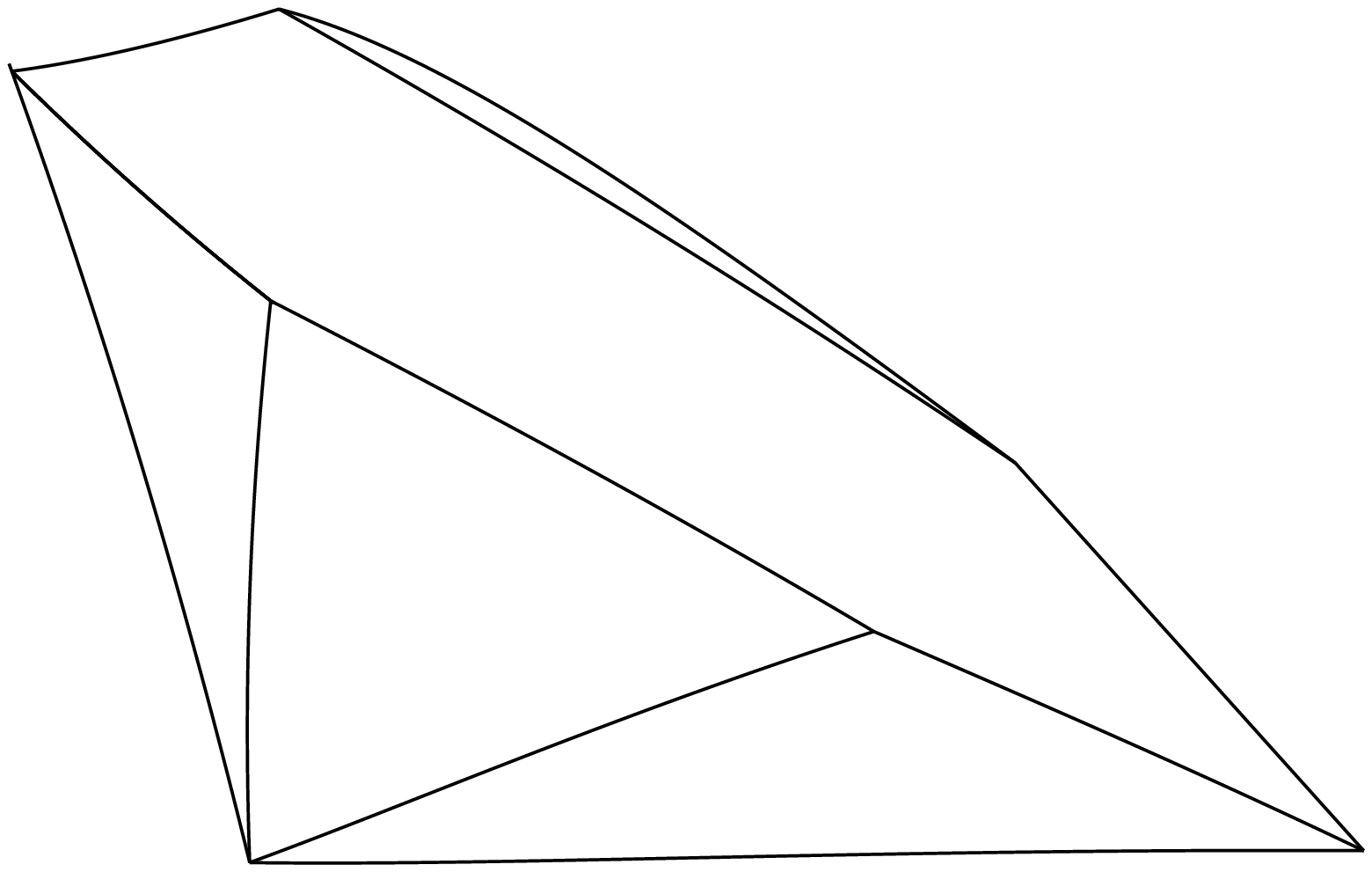, width=0.175\textwidth} }\hfill
  \subfigure[$123$]{
    \epsfig{figure=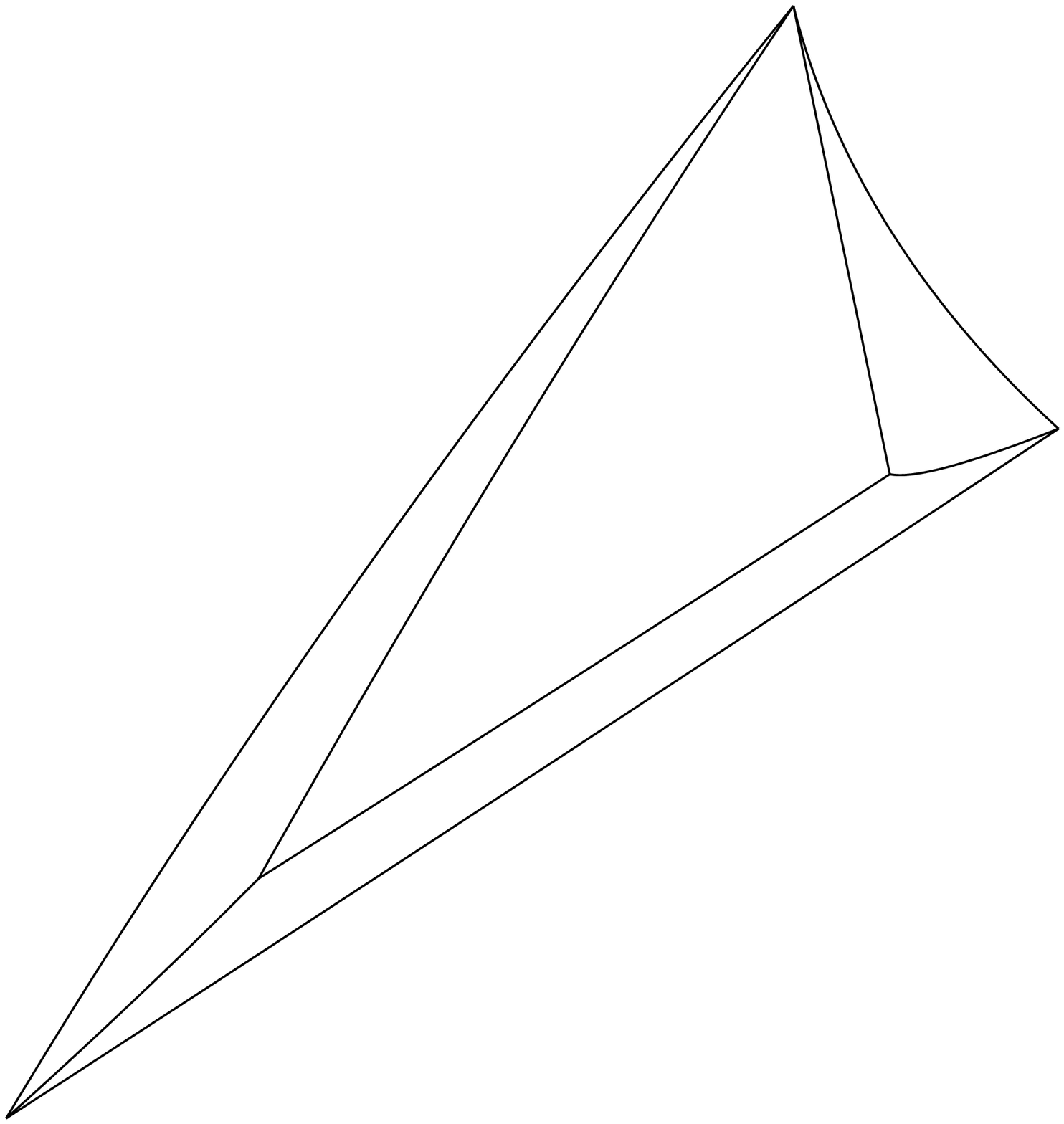, width=0.125\textwidth} }\hfill\mbox{}

  \hfill\subfigure[$1212$]{
    \epsfig{figure=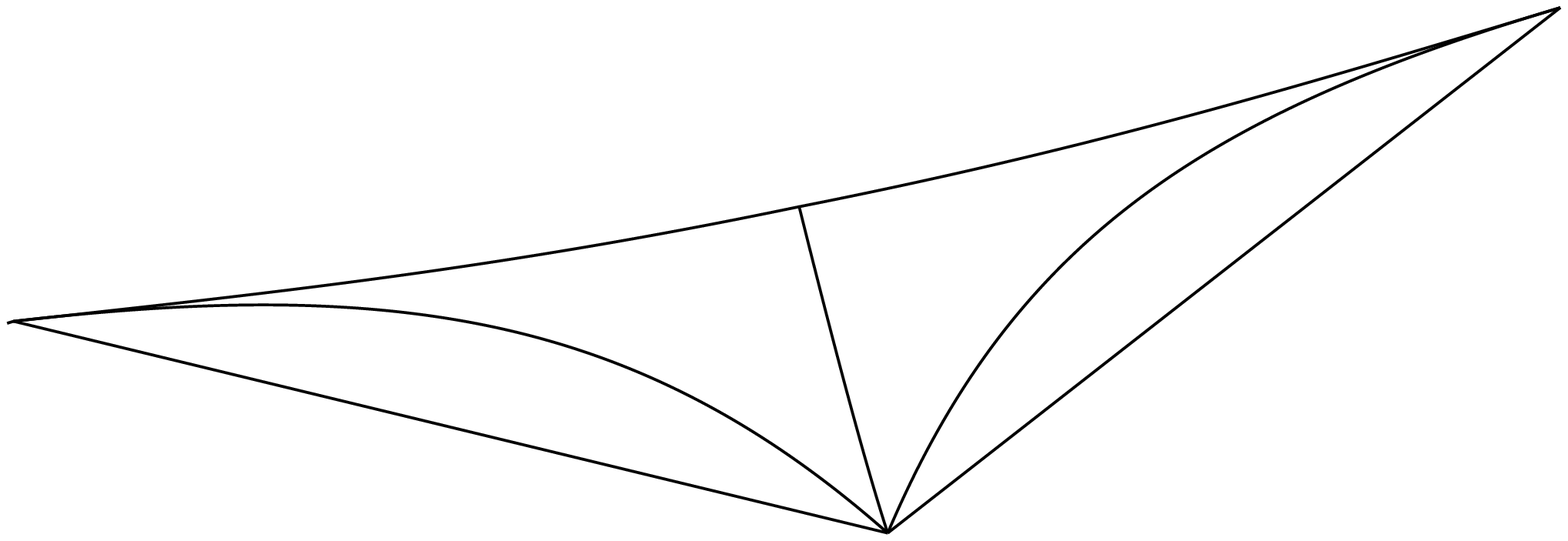, width=0.2\textwidth} }\hfill
  \subfigure[$12\overline{1}31$]{
    \epsfig{figure=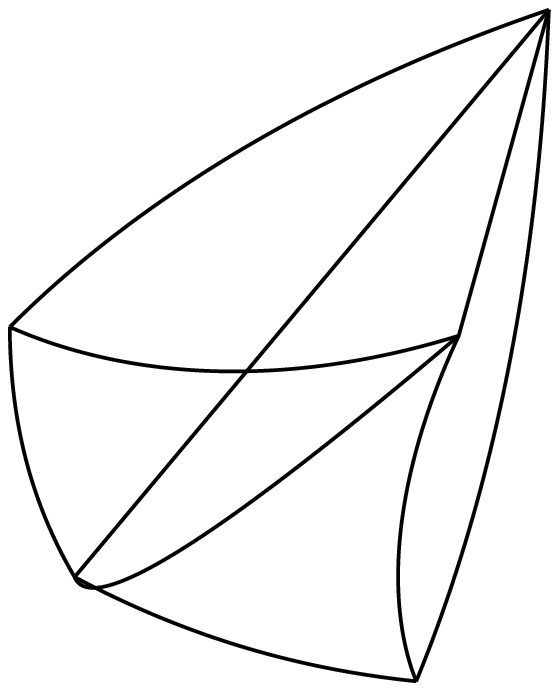, width=0.125\textwidth} }\hfill\mbox{}
  \caption{Faces of the Dirichlet domain for
    $\Gamma(2\pi/3,\overline{\sigma}_4)$, drawn in spinal coordinates.} \label{fig:pics3s4c}
\end{figure}

Pictures of the (isometry classes of) 3-faces of the Dirichlet domain
for $\Gamma(2\pi/3,\overline{\sigma}_4)$ are given in
Figure~\ref{fig:pics3s4c}. We chose to display the faces for that
specific group because its combinatorics are particularly simple among
all sporadic groups (Dirichlet domains for sporadic lattices can have
about a hundred faces).

In case of Behaviour A, the program provides a list of faces for the
polyhedron $F_W$, and checks whether it has side-pairings in the sense
of the Poincar\'e polyhedron theorem (once again, we choose to check
this only numerically). There is a minor issue of ambiguity between
the side pairings, due to the fact that most groups
$\Gamma(\frac{2\pi}{p},\tau)$ actually contain $J$, which means that
the centre of the Dirichlet domain has non-trivial
stabiliser. Possibly after adjusting the side-pairings by
pre-composing them with $J$ or $J^{-1}$, all the groups exhibiting
Behaviour~A turn out to have side-pairings (or at least they appear
to, numerically). Another way to take care of the issue of non-trivial
stabiliser for the centre of the Dirichlet domain is of course simply
to change the centre (within reasonably small distance to the centre
of mass of the mirrors, since we want the side-pairings obtained from
the Dirichlet domain to be related in simple terms to the original
generating reflections).

In either case, either after adjusting the side-pairings by elements
of the stabiliser, or after changing the centre, we are in a position
to check the cycle conditions of the Poincar\'e polyhedron theorem.
The general philosophy that grew out of~\cite{De1} (see
also~\cite{M1}, or even~\cite{Pic}) is that the only cycle conditions
that need to be checked are those for complex totally geodesic
$2$-faces, where the cycle transformations are simply complex
reflections. Our program goes through all these complex $2$-faces, and
computes the rotation angle of the cycle transformations (as well as
the total angle inside the polyhedron along the cycle).

Table~\ref{tab-poincare} gives the list of sporadic groups that
exhibit Behaviour~A and whose complex cycles rotate by an integer part
of $2\pi$ (for $\tau=\overline{\sigma}_4$, $p=8$, one needs to use a
centre for the Dirichlet domain other than the centre of mass of the
mirrors of the three reflections).
\begin{table}[htbp]
\centering
\begin{tabular}{|c|c|}\hline
$\tau$ & $p$\\
\hline
$\sigma_1$ & $p=3,4,6$\\ 
$\overline{\sigma}_4$ & $p=3,4,5,6,8,12$\\
$\sigma_5$ & $p=3,4$\\\hline
\end{tabular}
\caption{Sporadic groups with $3\leqslant p\leqslant 24$ whose Dirichlet domain satisfies the hypotheses of 
  the Poincar\'e polyhedron theorem, at least numerically.} \label{tab-poincare}
\end{table}

For groups that exhibit behaviour~A but whose cycle transformations
rotate by angles that are not integer parts of $2\pi$, all one can
quickly say is that the G-closed polyhedron cannot be a fundamental
domain for their action (even modulo the stabiliser of $p_0$), but the
group may still be a lattice. This issue is related to the question of
whether the integrality condition of~\cite{DM} is close to being
necessary and sufficient for the corresponding reflection group to be
a lattice (see the analysis in~\cite{M3}).

There is a natural refinement of the procedure described in
Section~\ref{sec-gprocedure} to handle this case. Suppose a given
cycle transformation $g$ rotates by an angle $\alpha$, and
$2\pi/\alpha$ is not an integer. If that number is not rational, the
group is not discrete (the irrationality can of course be difficult to
actually prove). If $\alpha=2\pi m/n$ for $m,n\in\Z$, then some power
$h=g^k$ rotates by an angle $2\pi/n$, and it is natural to replace the
G-closed set of group elements $W$ by 
\begin{equation}\label{eq:refinement}
W\cup h W h^{-1}.
\end{equation}
One then starts over with the G-procedure as described in
Section~\ref{sec-gprocedure}, starting from $W_0=W\cup h W h^{-1}$.

The groups with problematic rotation angles are
$$
\Gamma\left(\frac{2\pi}{5},\sigma_1\right),
\Gamma\left(\frac{2\pi}{5},{\sigma}_5\right),
$$
and all groups with $\tau=\overline{\sigma}_4$, $p\neq
3,4,5,6,8,12$. The ones with $\tau=\overline{\sigma}_4$ are known to
be non-discrete, see Theorem~\ref{thm:nondiscrete}. The groups
$\Gamma(\frac{2\pi}{5},\sigma_1)$ and
$\Gamma(\frac{2\pi}{5},{\sigma}_5)$ do not seem to be
discrete. Indeed, their Giraud-closed sets have problematic rotation
angles, see Table~\ref{tab-problematic}. In both cases, after
implementing the refinement of~\eqref{eq:refinement}, the G-procedure
exhibits Behaviour B.

\begin{table}[htbp] 
\centering
\begin{tabular}{r|c|c}
  Group & cycle transformation & angle\\ \hline
  $\Gamma(\frac{2\pi}{5},\sigma_1)$ & $(R_1R_2)^2$ & $4\pi/5$\\
  $\Gamma(\frac{2\pi}{5},\sigma_5)$ & $((R_1J)^5R_2^{-1})^2$ & $4\pi/15$
\end{tabular}
\caption{Some problematic rotation angles in Giraud-closed polyhedra.}\label{tab-problematic}
\end{table}

\section{Group presentations}

From the geometry of the Dirichlet domains for sporadic lattices, one
can infer explicit group presentations. Indeed one knows that the
side-pairings generate the group, and the relations are normally
generated by the cycle transformations, see~\cite{EP} for instance.

Given that there are many faces, it is of course quite prohibitive to
write down such a presentation by hand. It is reasonably easy however
to have a computer do this. Our program produces files that can be
passed to GAP in order to simplify the presentations (it is quite
painful, even though not impossible, to do these simplifications by
hand). It turns out that the presentations coming from the Dirichlet
domains can all be reduced to quite a simple form (see
Table~\ref{tab:presentations}).

Note that the results of this section are just as conjectural as the
statement of Conjecture~\ref{conj:main}, since they depend on the
accuracy of the combinatorics of the Dirichlet domains.

\begin{table}[htbp]
$$
\begin{array}{lc}
\Gamma(\frac{2\pi}{3},\sigma_1): 
& J=12312312=23123123=31231231; \\
& 1^3=Id;\quad (123)^8=Id; \\
& (12)^3=(21)^3;\quad [1(23\overline{2})]^2=[(23\overline{2})1]^2;\quad
1(232\overline{3}\overline{2})1
=(232\overline{3}\overline{2})1(232\overline{3}\overline{2}).\\
& \\
\Gamma(\frac{2\pi}{4},\sigma_1): 
& J=12312312=23123123=31231231; \\
& 1^4=Id;\quad (123)^8=Id; \quad (12)^{12}; \\
& (12)^3=(21)^3;\quad [1(23\overline{2})]^2=[(23\overline{2})1]^2;\quad
1(232\overline{3}\overline{2})1
=(232\overline{3}\overline{2})1(232\overline{3}\overline{2}). \\
& \\
\Gamma(\frac{2\pi}{6},\sigma_1):
& J=12312312=23123123=31231231; \\
& 1^6=Id;\quad (123)^8=Id; \quad (12)^{6};\quad [1(23\overline{2})]^{12}=Id; \\
& (12)^3=(21)^3;\quad [1(23\overline{2})]^2=[(23\overline{2})1]^2;\quad
1(232\overline{3}\overline{2})1
=(232\overline{3}\overline{2})1(232\overline{3}\overline{2}); \\
& \\
\Gamma(\frac{2\pi}{3},\overline{\sigma}_4):
& J^{-1}=1231231=2312312=3123123;\\
& 1^3=Id;\quad (123)^7=Id;  \\
& (12)^2=(21)^2;\\
& \\
\Gamma(\frac{2\pi}{4},\overline{\sigma}_4):
& J^{-1}=1231231=2312312=3123123;\\
& 1^4=Id;\quad (123)^7=Id;  \\
& (12)^2=(21)^2; \\
& \\
\Gamma(\frac{2\pi}{5},\overline{\sigma}_4):
& J^{-1}=1231231=2312312=3123123;\\
& 1^5=Id;\quad (123)^7=Id;\quad (12)^{20};  \\
& (12)^2=(21)^2; \\
& \\
\Gamma(\frac{2\pi}{6},\overline{\sigma}_4):
& J^{-1}=1231231=2312312=3123123;\\
& 1^6=Id;\quad (123)^7=Id;\quad (12)^{12};  \\
& (12)^2=(21)^2; \\
& \\
\Gamma(\frac{2\pi}{8},\overline{\sigma}_4):
& J^{-1}=1231231=2312312=3123123;\\
& 1^8=Id;\quad (123)^7=Id;\quad (12)^{8};\quad [1(23\overline{2})]^{24};  \\
& (12)^2=(21)^2; \\
& \\
\Gamma(\frac{2\pi}{12},\overline{\sigma}_4):
& J^{-1}=1231231=2312312=3123123;\\
& 1^{12}=Id;\quad (123)^7=Id;\quad (12)^{6};\quad [1(23\overline{2})]^{12};  \\
& (12)^2=(21)^2; \\
& \\
\Gamma(\frac{2\pi}{3},\sigma_5):
& J^3=Id;\quad J1J^{-1}=2;\quad J2J^{-1}=3;\quad J3J^{-1}=1;\\
& 1^3=Id;\quad (123)^{10};\\
& (12)^2=(21)^2;\quad 1(23\overline{2})1(23\overline{2})1
=(23\overline{2})1(23\overline{2})1(23\overline{2}). \\
& \\
\Gamma(\frac{2\pi}{4},\sigma_5):
& J^3=Id;\quad J1J^{-1}=2;\quad J2J^{-1}=3;\quad J3J^{-1}=1;\\
& 1^4=Id;\quad (123)^{10};\quad 
(1\overline{3}23\overline{1}23\overline{2})^{12};  \\
& (12)^2=(21)^2;\quad 1(23\overline{2})1(23\overline{2})1
=(23\overline{2})1(23\overline{2})1(23\overline{2}). \\
& \\
\end{array}
$$
\caption{Conjectural presentations for the groups that appear in
  Conjecture~\ref{conj:main}. The groups with $\tau=\sigma_1$,
  $\overline{\sigma}_4$ are generated by $R_1$, $R_2$ and $R_3$, that
  is $J$ can be expressed as a product of the $R_j$'s. For
  $\tau=\sigma_5$ this is not the case, and $\langle
  R_1,R_2,R_3\rangle$ has index $3$ in $\langle J, R_1\rangle$.}\label{tab:presentations}
\end{table}

\section{Description of the cusps of the non-compact examples}
\label{sec:cusps}

The geometry of the Dirichlet domains for sporadic lattices gives
information about the isotropy groups of any vertex. Rather than
giving a whole list, we gather information about the cusps in the
Dirichlet domain and in $M=\Gamma\setminus {\rm H}^2_\C$, by giving the
number of cusps, as well as generators and relations for their
stabilisers (see Table~\ref{tab:cusps}).  

Once again, the results of this section are conjectural (they depend on the
accuracy of the combinatorics of the Dirichlet domains).

\begin{table}[htbp]
\centering
\renewcommand{\arraystretch}{1.2}
\begin{tabular}{|c|c|c|c|l|ll|}\hline
 $p$ & $\tau$ & \# cusps & \# cusps in $M$ & Generators 
& Relations & \\ 
\hline
$3$ & $\sigma_1$ & 3 & 1 & $1$, $2$ & $1^3=2^3=Id$, & $(12)^3=(21)^3$ \\
\hline
$4$ & $\sigma_1$ & 6 & 1 & $1$, $23\overline{2}$ & 
$1^4=(23\overline{2})^4=Id$, & $[1(23\overline{2})]^2=[(23\overline{2})1]^2$ \\
\hline
$6$ & $\sigma_1$ & 6 & 2 & $1,\ 232\overline{3}\,\overline{2}$ & 
$1^6=(232\overline{3}\overline{2})^6=Id$, &
$1(232\overline{3}\overline{2})1
=(232\overline{3}\overline{2})1(232\overline{3}\overline{2})$ \\
    &            &   &   & $1,\ \overline{3}\,\overline{2}323$ & 
$1^6=(\overline{3}\overline{2}323)^6=Id$, &
$1(\overline{3}\overline{2}323)1
=(\overline{3}\overline{2}323)1(\overline{3}\overline{2}323)$ \\
\hline
$4$ & $\overline{\sigma}_4$ & 3 & 1 & $1$, $2$ & 
$1^4=2^4=Id$, & $(12)^2=(21)^2$ \\
\hline
$6$ & $\overline{\sigma}_4$ & 6 & 1 & $1$, $23\overline{2}$ & 
$1^6=(23\overline{2})^6=Id$, & 
$1(23\overline{2})1=(23\overline{2})1(23\overline{2})$\\
\hline
$3$ & $\sigma_5$ & 3 & 1 & $23\overline{2}$,\ $(1J)^5$ & 
$(23\overline{2})^3=[(IJ)^5]^6=Id$, & 
$[(23\overline{2})(1J)^{-5}]^2=[(1J)^{-5}(23\overline{2})]^2$ \\
\hline
$4$ & $\sigma_5$ & 3 & 1 & $1$, $2$ & 
$1^4=2^4=Id$, & $(12)^2=(21)^2$ \\
\hline
\end{tabular}
\caption{Conjectural list of cusps for the non-cocompact examples from
  Conjecture~\ref{conj:main}. All of the relations follow
directly from the conjectural presentations given in
Table~\ref{tab:presentations}; some follow directly but others with
slightly more work.}\label{tab:cusps}
\end{table}

\section{Non-discreteness results} \label{sec-nondiscrete}

In this section we prove some restrictions on the parameters for the
group $\Gamma(2\pi/p,\tau)$ to be discrete, aiming to show the
optimality of the statement of Conjecture~\ref{conj:main}. More
specifically, we will show the following.
\begin{thm} \label{thm:nondiscrete} Only finitely many of the sporadic triangle groups are discrete. More precisely:
\begin{itemize}
\item For $p \geqslant 7$, $\Gamma(\frac{2\pi}{p},\sigma_1)$ is not discrete.
\item For $p=3,5,6,7$, $\Gamma(\frac{2\pi}{p},\overline{\sigma}_1)$ is not discrete. 
\item For $p \geqslant 6$, $\Gamma(\frac{2\pi}{p},\sigma_2)$ is not discrete.
\item For $6 \leqslant p \leqslant 19$, $\Gamma(\frac{2\pi}{p},\overline{\sigma}_2)$ is not discrete.
\item For $p=4,5,6$, $\Gamma(\frac{2\pi}{p},\sigma_4)$ is not discrete.
\item For $p \neq 2,3,4,5,6,8,12$,  $\Gamma(\frac{2\pi}{p},\overline{\sigma}_4)$ is not discrete.
\item For $p \neq 2,3,4,5,6,8,12$, $\Gamma(\frac{2\pi}{p},\sigma_5)$ is not discrete.
\item $\Gamma(\frac{2\pi}{4},\overline{\sigma}_5)$ is not discrete.
\item $\Gamma(\frac{2\pi}{5},\sigma_6)$ and $\Gamma(\frac{2\pi}{5},\overline{\sigma}_6)$ are not discrete.
\item For $p \neq 2,3,4,7,14$, $\Gamma(\frac{2\pi}{p},\sigma_7)$ is not discrete.
\end{itemize}
\end{thm}

The proofs are slightly different for each part of the statement, as
detailed in Table~\ref{tab:nondisc}. Since all of them are based
either on Knapp's theorem or on J\o rgensen's inequality, we shall
briefly review these results in section~\ref{sec:review}.

\begin{table}[htbp]
$$
\begin{array}{|l|l|l|l|}
\hline
\tau & p \hbox{ with }H_\tau \hbox{ hyperbolic} 
& p \hbox{ where non-discrete} & \hbox{Result used}  \\
\hline
\sigma_1 & [3,\infty) & 7, 8, [10,\infty) & \hbox{Proposition \ref{prop:sig1}} \\
& & 9 & \hbox{Proposition \ref{prop:sig1'}} \\
\overline{\sigma}_1 & [3,7] & 5,7 & \hbox{Proposition \ref{prop:sig1c}} \\
& & 3,6,7 & \hbox{Proposition \ref{prop:sig1c'}} \\
\sigma_2 & [3,\infty) & [6,9], [11,\infty) & \hbox{Proposition \ref{prop:sig2}} \\
& & 10 & \hbox{Proposition \ref{prop:sig2'}} \\
\overline{\sigma}_2 & [3,19] & [6,9], [11,19] & \hbox{Proposition \ref{prop:sig2c}} \\
& & 10 & \hbox{Proposition \ref{prop:sig2'}} \\
\sigma_3 & [3,\infty) & [3,\infty) & \hbox{Proposition 4.5 of \cite{vol2}} \\
\overline{\sigma}_3 & [3,6] & [3,6] & \hbox{Proposition 4.5 of \cite{vol2}} \\
\sigma_4 & [4,6] & [4,6] & \hbox{Proposition \ref{prop:sig4}} \\
\overline{\sigma}_4 & [3,\infty) & 7, [9, 11], [13,\infty)
& \hbox{Proposition \ref{prop:sig4c}} \\
\sigma_5 & [2,\infty) & 7, [9,11],[13,\infty & \hbox{Proposition \ref{prop:sig5}} \\
\overline{\sigma}_5 & \{2,\,4\} & 4 & \hbox{Proposition \ref{prop:sig4}} \\
\sigma_6 & [3,\infty) & 3,4,[6\infty) & \hbox{Proposition 4.5 of \cite{vol2}} \\
& & 5 & \hbox{Proposition \ref{prop:sig4}} \\
\overline{\sigma}_6 & [3,29] & 3,4,[6\infty) & \hbox{Proposition 4.5 of \cite{vol2}} \\
& & 5 & \hbox{Proposition \ref{prop:sig4}} \\
\sigma_7 & [2,\infty) & 5,6,[8,13],[15,\infty) &  \hbox{Proposition \ref{prop:sig7}} \\
\overline{\sigma}_7 & \{2\} & & \\
\sigma_8 & [4,41] & [4,41] & \hbox{Corollary 4.2 of \cite{vol2}} \\
\overline{\sigma}_8 & [4,\infty) & [4,\infty) & \hbox{Corollary 4.2 of \cite{vol2}} \\
\sigma_9 & [3,\infty) & [3,\infty) & \hbox{Corollary 4.2 of \cite{vol2}} \\
\overline{\sigma}_9 & [4,8] & [4,8] & \hbox{Corollary 4.2 of \cite{vol2}} \\
\hline
\end{array}
$$
\caption{Values of the parameter where Knapp or J\o rgensen show
  non-discreteness. The second column gives the values of $p$ for
  which the Hermitian form $H_\tau$ has signature $(2,1)$ (taken from
  \cite{vol2}). The third and fourth columns give values of $p$ for
  which a well chosen subgroup fails the Knapp test or the J\o rgensen
  test (and hence the group is not discrete). If this was done in
  \cite{vol2} we give the reference. For some values of $\tau$ we
  apply Knapp and J\o rgensen to two different complex reflections in
  the group (in which case the results are listed on two separate
  lines).  }
\label{tab:nondisc}
\end{table}

\subsection{Knapp, J\o rgensen and Shimizu} \label{sec:review}

Knapp's theorem gives a necessary and sufficient condition for a
two-generator subgroup of ${\rm PU}(1,1)$ to be discrete, assuming
both generators as well as their product are elliptic. The reference
for Knapp's theorem is~\cite{Kna}, see also~\cite{KS}.  The full list
of possible rotation angles for $A$, $B$ and $AB$ will not be needed
here.  In fact we shall only use the following special case of Knapp's
theorem, that applies to isosceles triangles.

\begin{thm}\label{lem:knapp} (Knapp) Consider a triangle in ${\rm H}^2_\R$ with
  angles $\alpha,\alpha, \beta$, and let $\Delta$ be the group
  generated by the reflections in its sides. If $\Delta$ is discrete
  then one of the following holds:
\begin{itemize}
\item $\alpha=\frac{\pi}{q}$ and ($\beta=\frac{2\pi}{r}$ or
  $\frac{4\pi}{q}$) with $q,r \in \N^*$
\item $\alpha=\frac{2\pi}{r}$ and $\beta=\frac{2\pi}{r}$ with $r \in \N^*$.
\end{itemize}
\end{thm}

\begin{rmk}
In a few cases, we also use the spherical version of Knapp's
theorem, which is a result of Schwarz (see~\cite{vol2}).  
\end{rmk}

\begin{figure}[htbp]
  \centering
  \epsfig{figure=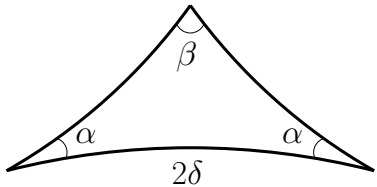, width=5cm}
  \caption{We shall apply Knapp's theorem in the special case of isosceles
triangles, see formula~(\ref{eq:knappthm}) for the relationship
between angles and distances.}
\label{fig:knapp}
\end{figure}

Basic hyperbolic trigonometry gives a relationship between the angles
and the length of the base of the triangle (see
Figure~\ref{fig:knapp}). Indeed, if the length of the base is
$2\delta$, then
\begin{equation}\label{eq:knappthm}
    \cosh \delta\ \sin \alpha = \cos\frac{\beta}{2}
\end{equation}
This gives a practical computational way to check whether the
conditions of Knapp's theorem hold.

Note also that the statement of Knapp's theorem implies that if
$\alpha=\pi/q$ for $q\in\N^*$, and if the angle $\beta$ is larger than
$2\pi/3$, then the group cannot be discrete. In view of
formula~\eqref{eq:knappthm}, the latter statement is the same 
as one would obtain from
J\o rgensen's inequality (see~\cite{JKP}):

\begin{thm} \label{thm:jorg} (Jiang-Kamiya-Parker) Let $A$ be a complex
  reflection through angle $2\alpha=\frac{2\pi}{q}$ with $q\in\N^*$, with
  mirror the complex line $L_A$.  Let $B \in {\rm PU}(2,1)$ be such that $B(L_A)$
  and $L_A$ are ultraparallel, and denote their distance by $2\delta$. If 
\begin{equation}\label{eq:jorgthm}
    |\cosh \delta\ \sin \alpha|< \frac{1}{2}
    \end{equation}
  then $\langle A, B \rangle$ is non-discrete.
\end{thm}

In certain cases we need to deal with groups generated by vertical
Heisenberg translations (see definition in section~\ref{sec-HnC}). In
this case we need results that generalise the above version of J\o
rgensen's inequality and Knapp's theorem.  These results are complex
hyperbolic versions of Shimizu's lemma, Proposition 5.2 of \cite{Par1}
and a lemma of Beardon, Theorem 3.1 of \cite{Par2}. We combine them in
the following statement which is equivalent to the statements given in
\cite{Par1} and \cite{Par2}.

\begin{thm}\label{thm:shim} (Parker)
Let $A\in{\rm SU}(2,1)$ be a parabolic map conjugate to a vertical Heisenberg 
translation with fixed point $z_A$. Let $B\in{\rm SU}(2,1)$ be a map not
fixing $z_A$. If $\langle A, B\rangle$ is discrete then
either ${\rm tr}(ABAB^{-1})\leqslant -1$ or ${\rm tr}(ABAB^{-1})=3-4\cos^2(\pi/r)$
for some $r\in\N$ with $r\geqslant 3$. In particular, if
\begin{equation}\label{eq:shim}
2<{\rm tr}(ABAB^{-1})<3
\end{equation}
then  $\langle A, B \rangle$ is non-discrete.
\end{thm}

\subsection{Using Knapp and J\o rgensen with powers of $R_1R_2$}

\subsubsection{The general set up}

Recall from~\cite{vol2} that for any sporadic value $\tau$, there is
a positive rational number $r/s$ so that
\begin{equation}\label{eq:taunormsquared}
|\tau|^2=2+2\cos(r \pi/s),
\end{equation}
which corresponds to the fact that $R_1R_2$ should have finite order. The
values of these $r$ and $s$ are clearly the same for $\sigma_j$ and
$\overline{\sigma}_j$, and are given by
\begin{equation}\label{eq:rands}
\begin{array}{|r|ccccccccc|}\hline
\tau &  \sigma_1 & \sigma_2 & \sigma_3 & \sigma_4 & \sigma_5 & \sigma_6 & \sigma_7 & \sigma_8 & \sigma_9\\\hline
r/s & 1/3 & 1/5 & 3/5 & 1/2 & 1/2 & 1/2 & 1/7 & 5/7 & 3/7\\ \hline
\end{array}  
\end{equation}
Straightforward calculation shows
$$
R_1R_2 = \left[\begin{matrix} 
e^{2i\pi/3p}(1-|\tau^2|) & e^{4i\pi/3p}\tau & \tau^2-\overline{\tau} \\
-\overline{\tau} & e^{2i\pi/3p} & e^{-2i\pi/3p}\tau \\ 
0 & 0 & e^{-4i\pi/3p} \end{matrix}\right]
$$
which has eigenvalues $-e^{2i\pi/3p}e^{ri\pi/s}$,
$-e^{2i\pi/3p}e^{-ri\pi/s}$, $e^{-4i\pi/3p}$.
Therefore $(R_1R_2)^s$ has a repeated eigenvalue. An
$e^{-4i\pi/3p}$-eigenvector of $R_1R_2$ is given by
$$
{\bf p}_{12} = \left[\begin{matrix} 
e^{-2i\pi/3p}\tau^2+e^{4i\pi/3p}\overline{\tau}
-e^{-2i\pi/3p}\overline{\tau} \\
e^{2i\pi/3p}\overline{\tau}^2+e^{-4i\pi/3p}\tau
-e^{2i\pi/3p}\tau \\
2\cos(2\pi/p)+2\cos(r\pi/s)
\end{matrix}\right].
$$
For most values of $p$ and $\tau$ 
this vector is negative, in which case its orthogonal complement (with
respect to $H_\tau$) gives a complex line in the ball. Hence (in most
cases) it is a complex reflection, and one checks easily that it
commutes with both $R_1$ and $R_2$.

Likewise, for most values of $p$, $\tau$, $(R_2R_3)^s$ is a complex
reflection that commutes with $R_2$ and $R_3$, and it fixes a complex
line whose polar vector is ${\bf p}_{23}=J({\bf p}_{12})$.
If the distance between these two lines is $2\delta_p$ then from Lemma~\ref{lem:twolines}:
\begin{equation}\label{eq:distance}
\cosh^2(\delta_p) =
\frac{\langle{\bf p}_{12},{\bf p}_{23}\rangle
\langle{\bf p}_{23},{\bf p}_{12}\rangle} 
{\langle{\bf p}_{12},{\bf p}_{12}\rangle 
\langle{\bf p}_{23},{\bf p}_{23}\rangle}
= \frac{\bigl|\overline{\tau}^2+e^{-2i\pi/p}\tau-\tau\bigr|^2}
{\bigl(2\cos(2\pi/p)+2\cos(r\pi/s)\bigr)^2}.
\end{equation}
The eigenvalues of $(R_1R_2)^s$ are $(-1)^{s+r}e^{2is\pi/3p}$,
$(-1)^{s-r}e^{2is\pi/3p}$, $e^{-4is\pi/3p}$. Therefore the rotation angle
of $(R_1R_2)^s$ is $(r+s)\pi+2s\pi/p$. This may or may not be of the
form $2\pi/c$. When it is not, we can find a positive integer $k$ so that
$(R_1R_2)^{sk}$ is a complex reflection whose angle has the form $2\pi/c$.
We define $2\alpha_p$ to be the smallest positive rotation angle among all
powers of $(R_1R_2)^s$.

Assuming that the parameter $\tau$ is fixed, the group
$\Gamma(2\pi/p,\tau)$ is indiscrete thanks to the J\o rgensen
inequality, for the values of $p$ satisfying:
\begin{equation}\label{eq:jorg}
\cosh \delta_p \sin\alpha_p 
=\frac{\bigl|\overline{\tau}^2+e^{-2i\pi/p}\tau-\tau\bigr|\sin\alpha_p}
{\bigl|2\cos(2\pi/p)+2\cos(r\pi/s)\bigr|}
< \frac{1}{2}.
\end{equation}
Likewise, in order to prove non-discreteness using Knapp's theorem we
seek values of $p$ for which
\begin{equation}\label{eq:knapp}
\cosh \delta_p \sin\alpha_p 
=\frac{\bigl|\overline{\tau}^2+e^{-2i\pi/p}\tau-\tau\bigr|\sin\alpha_p}
{\bigl|2\cos(2\pi/p)+2\cos(r\pi/s)\bigr|}
\neq \cos(\pi/q) \hbox{ or } \cos(2\alpha_p)
\end{equation}
for a natural number $q$.

Since $(R_1R_2)^s$ is a complex reflection that rotates through angle 
$(r+s)\pi+2s\pi/p$, we can apply the test of J\o rgensen's inequality simply
to $(R_1R_2)^s$. As 
$\bigl|\sin\bigl((r+s)\pi+2s\pi/p\bigr)\bigr|=\bigl|2\sin(2s\pi/p)\bigr|$
this involves finding values of $p$ for which 
\begin{equation}\label{eq-jorgold}
\Bigl|\cosh(\delta_p)\sin(2s\pi/p)\Bigr|
=\frac{\bigl|\overline{\tau}^2+e^{-2i\pi/p}\tau-\tau\bigr|\ 
|\sin(2s\pi/p)|}
{|2\cos(2\pi/p)+2\cos(r\pi/s)|}<\frac{1}{2}.
\end{equation}
For fixed $r$ and $s$, as $p$ tends to infinity, the left hand side
tends to zero. This shows at once that there can be only finitely many
discrete groups among all sporadic groups; the rest of the paper is
devoted to the proof of Theorem~\ref{thm:nondiscrete}, which is a vast
refinement of that statement.

In the next few sections, we shall apply Knapp or J\o rgensen to
various powers of $R_1R_2$ (other elements in the group as well) in
order to get the better non-discreteness results.

\subsubsection{Cases where $|\tau|^2=2$}

From~\eqref{eq:taunormsquared}, for any $\tau$ with $|\tau|^2=2$, we
have $r/s=1/2$ and so
$$
\cosh\delta_p=
\frac{\bigl|\overline{\tau}^2+e^{-2i\pi/p}\tau-\tau\bigr|}
{\bigl|2\cos(2\pi/p)\bigr|}.
$$
This happens for $\tau=\sigma_4$, $\overline{\sigma}_4$,  
$\sigma_5$, $\overline{\sigma}_5$,  
$\sigma_6$ or $\overline{\sigma}_6$. 
For all these values we have
$$
(R_1R_2)^2 = \left[\begin{matrix} 
-e^{4i\pi/3p} & 0 & 
e^{2i\pi/3p}\overline{\tau}+e^{-4i\pi/3p}\tau^2-e^{-4i\pi/3p}\overline{\tau} \\
0 & -e^{4i\pi/3p} & 
e^{-2i\pi/p}\tau + \overline{\tau}^2-\tau \\ 
0 & 0 & e^{-8i\pi/3p} \end{matrix}\right],
$$
which is a complex reflection commuting with both $R_1$ and $R_2$, and
whose rotation angle is $(p-4)\pi/p$. Note that $(p-4)\pi/p = 2\pi/c$
for some $c\in\Z\cup\{\infty\}$ if and only if $p$ and $c$ are as
given in the following table:
$$
\begin{array}{|l|rrrrrrr|}
\hline
p & 2 & 3 & 4 & 5 & 6 & 8 & 12 \\
c & -2 & -6 & \infty & 10 & 6 & 4 & 3 \\
\hline
\end{array}
$$
(When $p=4$, and hence $c=\infty$, we find that $(R_1R_2)^2$ is parabolic.)
For other values of $p$, by choosing an appropriate power $k$, we can
arrange that $(R_1R_2)^{2k}$ rotates by a smaller angle than
$(R_1R_2)^2$:

\begin{lem}\label{lem:alpha_p}
Let $(R_1R_2)^2$ be as above. There exists $k\in\Z$ so that 
$(R_1R_2)^{2k}$ has rotation angle $2\alpha_p$ where
$$
\alpha_p=\frac{\gcd(p-4,2p)\pi}{2p}.
$$
In particular
\begin{itemize}
\item If $p\equiv 1$ (mod 2), then $\gcd(p-4,2p)=1$,
and so $\alpha_p=\frac{\pi}{2p}$.
\item If $p\equiv 2$ (mod 4), then $\gcd(p-4,2p)=2$,
and so $\alpha_p=\frac{\pi}{p}$.
\item If $p\equiv 4$ (mod 8), then $\gcd(p-4,2p)=8$,
and so $\alpha_p=\frac{4\pi}{p}$. 
\item If $p\equiv 0$ (mod 8), then $\gcd(p-4,2p)=4$,
and so $\alpha_p=\frac{2\pi}{p}$. 
\end{itemize}
\end{lem}

\Pf
We want to find $k$ so that $k(p-4)\pi/p$ reduced modulo $2\pi$ is
``minimal''. More precisely, we write this as
$$
k(p-4)\pi/p-2\pi l=2\alpha_p
$$
for $k\in\N^*$, $l\in\N$, and we want to find $\alpha_p$ of the form
$\pi/c$ for some $c\in\N$. The optimal value of $k$ depends on
arithmetic properties of $p$. Let $d=\gcd(p-4,2p)$ then we can find
integers $k$ and $l$ so that $k(p-4)-l(2p)=d$. This means that
$k(p-4)\pi/p-2\pi l=d\pi/p$ and so $\alpha_p=d\pi/2p$.
This proves the first assertion.

If we write $(p-4)=ad$ and $2p=bd$ then, eliminating $p$, we have
$2ad+8=bd$ and so $d=1,\,2,\,4$ or $8$. It is easy to check which
values of $p$ correspond to which value of $d$.
\EPf

In the case where $c=\infty$ the map $(R_1R_2)^2$ is parabolic. 
Up to multiplying by a cube root of unity, we have
$$
{\rm tr}\Bigl((R_1R_2)^2J(R_1R_2)^2J^{-1}\Bigr)
={\rm tr}\Bigl((R_1R_2)^2(R_2R_3)^2\Bigr)
=3-\bigl|\overline{\tau}^2+e^{-2i\pi/p}\tau-\tau\bigr|^2.
$$
Thus applying Theorem \ref{thm:shim} with $A=(R_1R_2)^2$ and $J=B$ we see
can prove non-discreteness by showing that
\begin{equation}\label{eq:shim3}
\bigl|\overline{\tau}^2+e^{-2i\pi/p}\tau-\tau\bigr|<1
\quad \hbox{ or } \quad
\bigl|\overline{\tau}^2+e^{-2i\pi/p}\tau-\tau\bigr|\neq 2\cos(\pi/r)
\end{equation}
with $r$ a natural number at least $3$.

Checking \eqref{eq:jorg}, \eqref{eq:knapp} or \eqref{eq:shim3} is
best done by a computer. 

\begin{prop}\label{prop:sig4c}
Let $\tau=\overline{\sigma}_4=(-1-i\sqrt{7})/2$ and so $r/s=1/2$. Then:
\begin{itemize}
\item If $p$ is odd then \eqref{eq:jorg} holds for $p\geqslant 7$;
\item If $p\equiv 2$ (mod $4$) then \eqref{eq:jorg} holds for $p\geqslant 10$;
\item If $p\equiv 4$ (mod $8$) then \eqref{eq:knapp} holds for $p=20$ and
\eqref{eq:jorg} holds for $p\geqslant 28$;
\item If $p\equiv 0$ (mod $8$) then \eqref{eq:jorg} holds for $p\geqslant 16$.
\end{itemize}
Thus for all the values of $p$ given above $\langle (R_1R_2)^2,J\rangle$ 
and hence $\Gamma(\frac{2\pi}{p},\overline{\sigma}_4)$ is not discrete.
\end{prop}

\Pf
For the sake of concreteness, we list some of the values in the following table
$$
\begin{array}{|l|l|l|c|}
\hline
p & \alpha_p & \cosh(\delta_p)\sin(\alpha_p) & \\
\hline
7 & \pi/14 & 0.4257\dots & \eqref{eq:jorg} \\
9 & \pi/18 & 0.2650\dots & \eqref{eq:jorg} \\
\hline
10 & \pi/10 & 0.4423\dots & \eqref{eq:jorg} \\
14 & \pi/14 & 0.2774\dots & \eqref{eq:jorg} \\
\hline
20 & \pi/5 & 0.6748\dots & \eqref{eq:knapp} \\
28 & \pi/7 & 0.4754\dots & \eqref{eq:jorg} \\
\hline
16 & \pi/8 & 0.4601\dots & \eqref{eq:jorg} \\
24 & \pi/12 & 0.2889\dots & \eqref{eq:jorg} \\
\hline
\end{array}
$$
\EPf

\begin{prop}\label{prop:sig5}
Let $\tau=\sigma_5=e^{2i\pi/9}+e^{-i\pi/9}2\cos(2\pi/5)$ and so $r/s=1/2$.
Then:
\begin{itemize}
\item If $p$ is odd then~\eqref{eq:jorg} holds when $p\geqslant 7$; 
\item If $p\equiv 2$ (mod $4$) then~\eqref{eq:jorg} holds when $p\geqslant 10$; 
\item If $p\equiv 4$ (mod $8$) then~\eqref{eq:knapp} holds when $p=20$ 
and~\eqref{eq:jorg} holds when $p\geqslant 28$; 
\item If $p\equiv 0$ (mod $8$) then~\eqref{eq:jorg} holds when  $p\geqslant 16$. 
\end{itemize}
Thus for all the values of $p$ given above $\langle (R_1R_2)^2,J\rangle$
and hence $\Gamma(\frac{2\pi}{p},\sigma_5)$  is not discrete.
\end{prop}

\Pf
Some values are given in the following table:
$$
\begin{array}{|l|l|l|c|}
\hline
p & \alpha_p & \cosh(\delta_p)\sin(\alpha_p) & \\
\hline
7 & \pi/14 & 0.4977\dots & \eqref{eq:jorg} \\
9 & \pi/18 & 0.3011\dots & \eqref{eq:jorg} \\
\hline
10 & \pi/10 & 0.4974\dots & \eqref{eq:jorg} \\
14 & \pi/14 & 0.3032\dots & \eqref{eq:jorg} \\
\hline
20 & \pi/5 & 0.7202\dots & \eqref{eq:knapp} \\
28 & \pi/7 & 0.4988\dots & \eqref{eq:jorg} \\
\hline
16 & \pi/8 & 0.4980\dots & \eqref{eq:jorg} \\
24 & \pi/12 & 0.3053\dots & \eqref{eq:jorg} \\
\hline
\end{array}
$$
\EPf

Recall from \cite{vol2} that $\Gamma(\frac{2\pi}{p},\sigma_4)$ 
has signature $(2,1)$ exactly when $4 \leqslant p \leqslant 6$; that
 $\Gamma(\frac{2\pi}{p},\overline{\sigma}_5)$ has signature $(2,1)$
 exactly when $p=2$ or $4$, and that $\Gamma(\frac{2\pi}{p},\sigma_6)$
 and $\Gamma(\frac{2\pi}{p},\overline{\sigma}_6)$ are not discrete
 except possibly when $p=5$. Hence for each of these values of $\tau$ we only have
finitely many things to check. We gather these cases into a single result.

\begin{prop}\label{prop:sig4}
\begin{itemize}
\item If $\tau=\sigma_4=(-1+i\sqrt{7})/2$, and so $r/s=1/2$, 
and $p=4$ then~\eqref{eq:shim3} holds;
\item If $\tau=\sigma_4=(-1+i\sqrt{7})/2$, and so $r/s=1/2$,
and $p=5$ then~\eqref{eq:jorg} holds;
\item If $\tau=\sigma_4=(-1+i\sqrt{7})/2$, and so $r/s=1/2$, 
and $p=6$ then~\eqref{eq:knapp} holds;
\item If $\tau=\overline{\sigma}_5=e^{-2i\pi/9}+e^{i\pi/9}2\cos(2\pi/5)$, 
and so $r/s=1/2$,  and $p=4$ then~\eqref{eq:shim3} holds;
\item If $\tau=\sigma_6=e^{2i\pi/9}+e^{-i\pi/9}2\cos(4\pi/5)$, 
and so $r/s=1/2$, and $p=5$ then~\eqref{eq:knapp} holds;
\item If $\tau=\overline{\sigma}_6=e^{-2i\pi/9}+e^{i\pi/9}2\cos(4\pi/5)$, 
and so $r/s=1/2$, and $p=5$ then~\eqref{eq:knapp} holds.
\end{itemize}
Thus for these values of $\tau$ and $p$ then $\langle (R_1R_2)^2,J\rangle$
and hence $\Gamma(\frac{2\pi}{p},\tau)$  is not discrete.
\end{prop}

\Pf
Suppose $\tau=\sigma_4=(-1+i\sqrt{7})/2$. If $p=4$ we have
$$
\bigl|\overline{\tau}^2+e^{-2i\pi/p}\tau-\tau\bigr|
=\sqrt{3-\sqrt{7}}=0.595\ldots
$$
If $p=5$ then $\alpha_p=\pi/10$ and $\cosh(\delta_p)\sin(\alpha_p)=0.445\ldots$.
If $p=6$ then $\alpha_p=\pi/6$ and 
$\cosh(\delta_p)\sin(\alpha_p)=0.550\ldots\in\bigl(\cos(\pi/3),\cos(\pi/4)\bigr)$.

If $\tau=\overline{\sigma}_5=e^{-2i\pi/9}+e^{i\pi/9}2\cos(2\pi/5)$ and $p=4$ then
$$
\bigl|\overline{\tau}^2+e^{-2i\pi/p}\tau-\tau\bigr|
=\sqrt{\frac{7+\sqrt{5}-3\sqrt{3}-\sqrt{15}}{2}}
=0.289\ldots .
$$ 

If $\tau=\sigma_6=e^{2i\pi/9}+e^{-i\pi/9}2\cos(4\pi/5)$ and
$p=5$ then $\cosh(\delta_p)\sin(\alpha_p)=0.937\ldots\in\bigl(\cos(\pi/8),\cos(\pi/9)\bigr)$.

If $\tau=\overline{\sigma}_6=e^{-2i\pi/9}+e^{i\pi/9}2\cos(4\pi/5)$ and
$p=5$ then $\cosh(\delta_p)\sin(\alpha_p)=0.750\ldots\in\bigl(\cos(\pi/4),\cos(\pi/5)\bigr)$.
\EPf

\subsubsection{Cases where $|\tau|^2=3$}

We now consider the case $|\tau|^2=3$, which happens for
$\tau=\sigma_1$ or $\overline{\sigma}_1$. In this case $r/s=1/3$ and 
$$
(R_1R_2)^3=\left[\begin{matrix}
e^{2i\pi/p} & 0 & (e^{-2i\pi/p}-1)\bigl(e^{-2i\pi/3p}\tau^2
+e^{4i\pi/3p}\overline{\tau}-e^{-2i\pi/3p}\overline{\tau}\bigr) \\
0 & e^{2i\pi/p} & (e^{-2i\pi/p}-1)
\bigl(e^{2i\pi/3p}\tau^2+e^{-4i\pi/p}\tau-e^{2i\pi/3p}\tau\bigr) \\
0 & 0 & e^{-4i\pi/p} \end{matrix}\right].
$$
This is a complex reflection commuting with both $R_1$ and $R_2$, with
angle $6\pi/p$.
As above, we want to check whether~\eqref{eq:jorg} holds for
$\alpha_p$ the smallest possible rotation angle of powers of
$(R_1R_2)^3$.
\begin{itemize}
\item
If $p\equiv 1$ or $2$ (mod $3$), then we can find $k,l\in\N$ so that
$6k\pi/p-2\pi l=2\pi/p$. Hence $\alpha_p=\pi/p$.
\item If $3$ divides $p$ then $6\pi/p$ is already in the form $2\pi/c$, hence
$\alpha_p=3\pi/p$.
\end{itemize}

\begin{prop}\label{prop:sig1}
Let $\tau=\sigma_1=e^{i\pi/3}+e^{-i\pi/6}2\cos(\pi/4)$ and so $r/s=1/3$. Then
\begin{itemize}
\item If $p\equiv 1$ or $2$ (mod $3$) 
then~\eqref{eq:knapp} holds when $p=7$ and~\eqref{eq:jorg} holds when 
$p\geqslant 8$; 
\item If $p$ is divisible by $3$ then~\eqref{eq:knapp} holds when $p=12$, $15$
or $18$ and~\eqref{eq:jorg} holds when $p\geqslant 21$. 
\end{itemize}
Thus for all the values of $p$ given above $\langle (R_1R_2)^3,J\rangle$ 
and hence $\Gamma(\frac{2\pi}{p},\sigma_1)$ is not discrete.
\end{prop}

\Pf
Some values are given in the following table
$$
\begin{array}{|l|l|l|c|}
\hline
p & \alpha_p & \cosh(\delta_p)\sin(\alpha_p) & \\
\hline
7 & \pi/7 & 0.6510\dots & \eqref{eq:knapp} \\
8 & \pi/8 & 0.4969\dots & \eqref{eq:jorg} \\
\hline
12 & \pi/4 & 0.8134\dots & \eqref{eq:knapp} \\
15 & \pi/5 & 0.6510\dots & \eqref{eq:knapp} \\
18 & \pi/6 & 0.5416\dots & \eqref{eq:knapp} \\
21 & \pi/7 & 0.4631\dots & \eqref{eq:jorg} \\
\hline
\end{array}
$$
\EPf

From \cite{vol2} we know that if 
$\tau=\overline{\sigma}_1=e^{-i\pi/3}+e^{i\pi/6}2\cos(\pi/4)$ then the only values of 
$p$ that give signature $(2,1)$ are those with $3\leqslant  p\leqslant 7$. 

\begin{prop}\label{prop:sig1c}
Let $\tau=\overline{\sigma}_1=e^{-i\pi/3}+e^{i\pi/6}2\cos(\pi/4)$ and so $r/s=1/3$.
\begin{itemize}
\item If $p=5$ then~\eqref{eq:knapp} holds;
\item If $p=7$ then~\eqref{eq:jorg} holds.
\end{itemize}
Thus for $p=5$ or $7$ we see that $\langle (R_1R_2)^3,J\rangle$ 
and hence $\Gamma(\frac{2\pi}{p},\overline{\sigma}_1)$ is not discrete.
\end{prop}

\subsubsection{Cases where $|\tau|^2=2+2\cos(\pi/5)$}

This happens for $\tau=\sigma_2$ or $\overline{\sigma}_2$. In that
case $r/s=1/5$ and $(R_1R_2)^5$ is a complex reflection with eigenvalues
$e^{10i\pi/3p}$, $e^{10i\pi/3p}$, $e^{-20i\pi/3p}$, thus it has
rotation angle $10\pi/p$. 
\begin{itemize}
\item If $p$ is not divisible by $5$, then we can find $k,l\in\N$ such
  that $10k\pi/p-2\pi l=2\pi/p$. Hence $\alpha_p=\pi/p$.
\item If $p$ is divisible by $5$, then $10\pi/p$ is already in the
  form $2\pi/c$, hence $\alpha_p=5\pi/p$.
\end{itemize}

\begin{prop}\label{prop:sig2}
Let $\tau=\sigma_2=e^{i\pi/3}+e^{-i\pi/6}2\cos(\pi/5)$ and so $r/s=1/5$.
\begin{itemize}
\item If $p$ is not divisible by $5$ then~\eqref{eq:knapp} holds when $p=6$ or $7$
and~\eqref{eq:jorg} holds when $p\geqslant 8$;
\item If $p$ is divisible by $5$ then~\eqref{eq:knapp} holds when 
$15\leqslant p\leqslant 30$ and ~\eqref{eq:jorg} holds when $p\geqslant 35$.
\end{itemize}
Thus for for these values of $p$ we see that $\langle (R_1R_2)^5,J\rangle$,  
and hence, $\Gamma(\frac{2\pi}{p},\sigma_2)$ is not discrete.
\end{prop}

\Pf
Some values are given in the following table
$$
\begin{array}{|l|l|l|c|}
\hline
p & \alpha_p & \cosh(\delta_p)\sin(\alpha_p) & \\
\hline
6 & \pi/6 & 0.631\ldots & \eqref{eq:knapp} \\
7 & \pi/7 & 0.516\ldots & \eqref{eq:knapp} \\
8 & \pi/8 & 0.438\ldots & \eqref{eq:jorg} \\
\hline
15 & \pi/3 & 0.908\ldots & \eqref{eq:knapp} \\
20 & \pi/4 & 0.729\ldots & \eqref{eq:knapp} \\
25 & \pi/5 & 0.601\ldots & \eqref{eq:knapp} \\
30 & \pi/6 & 0.508\ldots & \eqref{eq:knapp} \\
35 & \pi/7 & 0.440\ldots & \eqref{eq:jorg} \\
\hline
\end{array}
$$
\EPf.

From \cite{vol2} we know that if 
$\tau=\overline{\sigma}_2=e^{-i\pi/3}+e^{i\pi/6}2\cos(\pi/5)$ then the Hermitian
form has signature $(2,1)$ only when $3\le p\le 19$.

\begin{prop}\label{prop:sig2c}
Let  $\tau=\overline{\sigma}_2=e^{-i\pi/3}+e^{i\pi/6}2\cos(\pi/5)$ and so $r/s=1/5$.
\begin{itemize}
\item If $p$ is not divisible by $5$ then~\eqref{eq:knapp} holds when $p=6$ 
and~\eqref{eq:jorg} holds when $7\leqslant p\leqslant 19$;
\item If $p$ is divisible by $5$ then~\eqref{eq:knapp} holds when 
$p=15$.
\end{itemize}
Thus for for these values of $p$ we see that $\langle (R_1R_2)^5,J\rangle$,  
and hence $\Gamma(\frac{2\pi}{p},\overline{\sigma}_2)$, is not discrete.
\end{prop}

\Pf
Some values are given below
$$
\begin{array}{|l|l|l|c|}
\hline
p & \alpha_p & \cosh(\delta_p)\sin(\alpha_p) & \\
\hline
6 & \pi/6 & 0.5660 & \eqref{eq:knapp} \\
7 & \pi/7 & 0.4713 & \eqref{eq:jorg} \\
\hline
15 & \pi/5 & 0.8718 & \eqref{eq:knapp} \\
\hline
\end{array}
$$
\EPf

\subsubsection{Cases where $|\tau|^2=2+2\cos(\pi/7)$}

This happens for $\tau=\sigma_7$ or $\overline{\sigma}_7$. In this case
$r/s=1/7$. The only group with $\tau=\overline{\sigma}_7$ and signature $(2,1)$
is $p=2$. This group is a relabelling of the group with $\tau=\sigma_7$
and $p=2$. It is discrete.
So for the remainder of this section we consider the case
when $\tau=\sigma_7=e^{2i\pi/9}+e^{-i\pi/9}2\cos(2\pi/7)$.

Then $(R_1R_2)^7$ is a complex reflection with eigenvalues
$e^{14i\pi/3p}$, $e^{14i\pi/3p}$, $e^{-28i\pi/3p}$; thus it has rotation angle
$14\pi/p$. 
\begin{itemize}
\item If $p$ is not divisible by $7$, then we can find $k,l\in\N$ so
  that $14k\pi/p-2\pi l=2\pi/p$. Hence $\alpha_p=\pi/p$.
\item If $p$ is divisible by $7$, then $14\pi/p$ is already in the
  form $2\pi/c$, hence $\alpha_p=7\pi/p$.
\end{itemize}

\begin{prop}\label{prop:sig7}
Let $\tau=\sigma_7=e^{2i\pi/9}+e^{-i\pi/9}2\cos(2\pi/7)$ and so $r/s=1/7$.
\begin{itemize}
\item If $p$ is not divisible by 7 then~\eqref{eq:knapp} holds for 
$p=5$ or $6$ and~\eqref{eq:jorg} holds when $p\geqslant 8$;
\item If $p$ is divisible by $7$ then~\eqref{eq:knapp} holds for 
$21\leqslant p\leqslant 42$ and~\eqref{eq:jorg} holds when $p\geqslant 49$;
with $p\geqslant 49$.
\end{itemize}
Thus for for these values of $p$ we see that $\langle (R_1R_2)^7,J\rangle$, 
and hence $\Gamma(\frac{2\pi}{p},\sigma_7)$, is not discrete.
\end{prop}

\Pf
Some values are given below
$$
\begin{array}{|l|l|l|c|}
\hline
p & \alpha_p & \cosh(\delta_p)\sin(\alpha_p) & \\
\hline
5 & \pi/5 & 0.929\ldots & \eqref{eq:knapp} \\
6 & \pi/6 & 0.702\ldots & \eqref{eq:knapp} \\
8 & \pi/8 & 0.476\ldots & \eqref{eq:jorg} \\
\hline
21 & \pi/3 & 0.921\ldots & \eqref{eq:knapp} \\
28 & \pi/4 & 0.739\ldots & \eqref{eq:knapp} \\
35 & \pi/5 & 0.608\ldots & \eqref{eq:knapp} \\
42 & \pi/6 & 0.514\ldots & \eqref{eq:knapp} \\
49 & \pi/7 & 0.444\ldots & \eqref{eq:jorg} \\
\hline
\end{array}
$$
\EPf

\subsection{Using Knapp and J\o rgensen with powers of
$R_1R_2R_3R_2^{-1}$}

\subsubsection{The general set up}

A starightforward calculation shows that:
$$
R_1R_2R_3R_2^{-1}=\left[\begin{matrix}
e^{2i\pi/3p}(1-|\tau^2-\overline{\tau}|^2) &
e^{-2i\pi/3p}(\tau-(\tau^2-\overline{\tau})\overline{\tau}) &
-\tau^2+(\tau^2-\overline{\tau})(|\tau|^2-e^{2i\pi/p}) \\
\tau(\tau-\overline{\tau}^2) &
e^{-4i\pi/3p}(1-|\tau|^2) & 
e^{-2i\pi/3p}\tau(|\tau|^2-1+e^{2i\pi/p}) \\
e^{-2i\pi/3p}(\tau-\overline{\tau}^2) &
-e^{-2i\pi/p}\overline{\tau} &
e^{2i\pi/3p}+e^{-4i\pi/3p}|\tau|^2
\end{matrix}\right],
$$
hence ${\rm  tr}(R_1R_2R_3R_2^{-1})=e^{2i\pi/3p}(2-|\tau^2-\overline{\tau}|^2)
+e^{-4i\pi/3p}$. An $e^{-4i\pi/3p}$ eigenvector of $R_1R_2R_3R_2^{-1}$
is given by
$$
{\bf p}_{1232}=\left[\begin{matrix}
e^{-4i\pi/3p}
\bigl(\tau(1-e^{2i\pi/p})-(\tau^2-\overline{\tau})\overline{\tau}\bigr) \\
|\tau|^2(1-e^{-2i\pi/p})-\overline{\tau}(\overline{\tau}^2-\tau)
-|1-e^{2i\pi/p}|^2 \\
e^{-2i\pi/p}
\bigl(\overline{\tau}(1-e^{-2i\pi/p})-(\overline{\tau}^2-\tau)\tau\bigr)
\end{matrix}\right].
$$
Suppose that $|\tau^2-\overline{\tau}|^2=2+2\cos(r'\pi/s')$. Then
$(R_1R_2R_3R_2^{-1})^s$ is a complex reflection. The values of $r'$ and $s'$ 
are clearly the same for $\sigma_j$ and $\overline{\sigma}_j$. 
They are (see \cite{vol2}):
\begin{equation}\label{eq:r'ands'}
\begin{array}{|r|cccccc|}\hline
\tau &  \sigma_1 & \sigma_2 & \sigma_3 & \sigma_4 & \sigma_5 & \sigma_6 \\\hline
r'/s' & 1/2 & 1/3 & 1/3 & 2/3 & 2/5 & 4/5 \\ \hline
\end{array}  
\end{equation}

Let $2\delta'_p$ denote
the distance from its mirror to the image of its mirror under $J$
(with polar vector ${\bf p}_{2313}=J({\bf p}_{1232})$). Then, from Lemma~\ref{lem:twolines}:
$$
\cosh(\delta'_p)=\frac{|\langle{\bf p}_{2313},{\bf p}_{1232}\rangle|}
{|\langle{\bf p}_{1232},{\bf p}_{1232}\rangle|}
=\frac{\bigl|(1-e^{2i\pi/p})\tau
+|\tau|^2(\overline{\tau}^2-2\tau)+e^{-2i\pi/p}\overline{\tau}^2\bigr|}
{2\cos(2\pi/p)+2\cos(r'\pi/s')}.
$$

Let $\alpha'_p$ be the smallest non-zero angle that a power of 
$(R_1R_2R_3R_2^{-1})^s$ rotates by. Let $\delta'_p$, $r'$ and $s'$
be as above. 
In order to prove non-discreteness using the J\o rgensen inequality,
we need to find values of $p$ such that
\begin{equation}\label{eq:jorg'}
\cosh \delta'_p \sin\alpha'_p 
=\frac{\bigl|(1-e^{2i\pi/p})\tau
+|\tau|^2(\overline{\tau}^2-2\tau)+e^{-2i\pi/p}\overline{\tau}^2\bigr|\sin\alpha'_p}
{\bigl|2\cos(2\pi/p)+2\cos(r'\pi/s')\bigr|}
< \frac{1}{2}.
\end{equation}
In order to prove non-discreteness using Knapp's theorem we must 
find values of $p$ for which
\begin{equation}\label{eq:knapp'}
\cosh \delta'_p \sin\alpha'_p 
=\frac{\bigl|(1-e^{2i\pi/p})\tau
+|\tau|^2(\overline{\tau}^2-2\tau)+e^{-2i\pi/p}\overline{\tau}^2\bigr|\sin\alpha'_p}
{\bigl|2\cos(2\pi/p)+2\cos(r'\pi/s')\bigr|}
\neq \cos(\pi/q) \hbox{ or } \cos(2\alpha_p)
\end{equation}
for a natural number $q$.

\subsection{When $|\tau^2-\overline{\tau}|^2=2$}

In this case $\tau=\sigma_1$ or $\overline{\sigma}_1$, and $r'/s'=1/2$. Moreover,
$(R_1R_2R_3R_2^{-1})^2$ is a complex reflection with
angle $(p-4)\pi/p$. So we proceed as in the section with $|\tau|^2=2$.
In particular, $\alpha'_p$ is given by Lemma \ref{lem:alpha_p}.

Using Proposition \ref{prop:sig1}, we already know that when $p=7$, $8$
or $p\ge 10$ then $\Gamma(\frac{2\pi}{p},\sigma_1)$ is not discrete.
Therefore, we restrict our attention to $p\le 9$.

\begin{prop}\label{prop:sig1'}
Let $\tau=\sigma_1=e^{i\pi/3}+e^{-i\pi/6}2\cos(\pi/4)$ and so $r'/s'=1/2$. Then
~\eqref{eq:knapp'} holds for $p=9$.
Thus 
$\langle (R_1R_2R_3R_1^{-1})^2,J\rangle$,
and hence also $\Gamma(\frac{2\pi}{p},\sigma_1)$, is not discrete.
\end{prop}

\Pf
In this case $\alpha'_p=\pi/18$ and 
$\cosh(\delta'_p)\sin(\alpha'_p)=0.686\ldots\in\bigl(\cos(\pi/3),\cos(\pi/4)\bigr)$. 
\EPf

For $\overline{\sigma}_1$, recall from \cite{vol2} that
$\Gamma(\frac{2\pi}{p},\overline{\sigma}_1)$ has signature $(2,1)$
exactly when $3 \leqslant p \leqslant 7$.

\begin{prop}\label{prop:sig1c'}
Let $\tau=\overline{\sigma}_1=e^{-i\pi/3}+e^{+i\pi/6}2\cos(\pi/4)$ and so
$r'/s'=1/2$. Then
\begin{itemize}
\item If $p=3$ or $p=6$ then~\eqref{eq:knapp'} holds; 
\item If $p=7$ then~\eqref{eq:jorg'} holds.
\end{itemize}
Thus for $p=3$, $6$ or $7$ the group $\langle (R_1R_2R_3R_1^{-1})^2,J\rangle$,
and hence also $\Gamma(\frac{2\pi}{p},\overline{\sigma}_1)$, is not discrete.
\end{prop}

\Pf
The values of $\cosh(\delta'_p)\sin(\alpha'_p)$ are:
$$
\begin{array}{|l|l|l|c|}
\hline
p & \alpha_p & \cosh(\delta'_p)\sin(\alpha'_p) & \\ 
\hline
3 & \pi/6 & 0.982\ldots & \eqref{eq:knapp'} \\
7 & \pi/14 & 0.269\ldots & \eqref{eq:jorg'} \\
\hline
6 & \pi/6 & 0.859\ldots & \eqref{eq:knapp'} \\
\hline
\end{array}
$$
\EPf

\subsection{Cases where $|\tau^2-\overline{\tau}|^2=3$}

We only consider the case $\tau=\sigma_2$ or $\overline{\sigma}_2$
(since $\sigma_3$ or $\overline{\sigma}_3$ were already handled
in~\cite{vol2}).  In this case $r'/s'=1/3$. 
Then $(R_1R_2R_3R_1^{-1})^3$ is a complex reflection
with angle $6\pi/p$. So we proceed as in the section with
$|\tau|^2=3$. Namely,
\begin{itemize}
\item if $p$ is not divisible by $3$, some power gives an angle $\alpha_p=\pi/p$;
\item if $p$ is divisible by $3$, some power gives $\alpha_p=3\pi/p$.
\end{itemize}
In order to use J\o rgensen, we check whether
$\cosh\delta'_p\sin\alpha_p<\frac{1}{2}$.

For $\tau=\sigma_2$, using Propositions \ref{prop:sig2} and \ref{prop:sig2c}
we only need to consider the cases where $p\le 5$ or $p=10$. This method yields
nothing new for $p\le 5$.

\begin{prop}\label{prop:sig2'}
Let $p=10$.
\begin{itemize}
\item If $\tau=\sigma_2=e^{i\pi/3}+e^{-i\pi/6}2\cos(\pi/5)$, and so
$r'/s'=1/3$, then~\eqref{eq:knapp'} holds; 
\item If $\tau=\overline{\sigma}_2=e^{-i\pi/3}+e^{i\pi/6}2\cos(\pi/5)$, and so
$r'/s'=1/3$, then~\eqref{eq:jorg'} holds. 
\end{itemize}
Thus for $p=10$ and $\tau=\sigma_2$ or $\overline{\sigma}_2$, the group 
$\langle (R_1R_2R_3R_2^{-1})^3,J\rangle$ is not discrete.
Hence $\Gamma(\frac{2\pi}{10},\sigma_2)$ and
$\Gamma(\frac{2\pi}{10},\overline{\sigma}_2)$ are not discrete.
\end{prop}

\Pf
When $p=10$ and $\tau=\sigma_2$ we have
$$
\cosh(\delta'_p)\sin(\alpha'_p)=0.6181\ldots\in\bigl(\cos(\pi/3),\cos(\pi/4)\bigr).
$$
When $p=10$ and $\tau=\overline{\sigma}_2$ we have
$\cosh(\delta'_p)\sin(\alpha'_p)=0.3871\ldots<1/2$.
\EPf

\raggedright 

\end{document}